\documentclass[11pt]{amsart}
\usepackage{latexsym,amssymb,amsthm}
\usepackage{amsmath}
\def\s{\sigma}
\def\t{\tau}
\def\upperp{ {}^{\perp}}
\def\tbase{{\rm Base}}

\def\cD{{\mathcal D}}

\def\cO{{\mathcal O}}
\def\CC{\mathbb C}

\def\ZZ{\mathbb Z}
\def\SS{\mathbb S}

\def\PP{\mathbb P}
\def\OO{\mathbb O}
\def\JO{{\mathcal J}_3(\OO)}

\def\QQ{\mathbb Q}
\def\FF{\mathbb F\mathbb F}
\def\fh{{\mathfrak h}}
\def\fs{{\mathfrak s}}
\def\fsl{{\mathfrak {sl}}}

\def\Spin{{Spin}}

\def\fe{{\mathfrak e}}
\def\fg{{\mathfrak g}}
\def\fn{{\mathfrak n}}
\def\fp{{\mathfrak p}}
\def\ft{{\mathfrak t}}
\def\fl{\mathfrak l}
\def\l{\lambda}
\def\a{\alpha}
\def\o{\omega}
\def\oo{\Omega}
\def\b{\beta}
\def\g{\gamma}

\def\d{\delta}
\def\m{\mu}
\def\e{\varepsilon}
\def\ot{{\mathord{\otimes }\,}}
\def\op{{\mathord{\oplus }\,}}
 
\def\lra{{\mathord{\;\longrightarrow\;}}} \def\ra{{\mathord{\;\rightarrow\;}}}
\def\we{\wedge}
\def\tr{{\rm trace}\,}

\def\CC{\mathbb C }\def\BC{\mathbb C }\def\BP{\mathbb P }
\def\BG{\mathbb G }
\def\BO{\mathbb O}\def\BS{\mathbb S }
\def\al#1{\alpha_{#1}}

\def\bcc#1{\mathbb C^{#1}}

\def\ff#1{\mathbb F\mathbb F^{#1}}
\def\fff#1#2{\mathbb F\mathbb F^{#1}_{#2}} 
\def\fg{\mathfrak g}
\def\fh{\mathfrak h}
\def\fp{\mathfrak p}
\def\ft{\mathfrak t}

\def\hd{, ... ,}
\def\inv{{}^{-1}}
\def\ii{| \FF^2_{X,x}|}

\def\La#1{\Lambda^{#1}}
\def\na{n+a}
\def\lo#1{\omega_{#1}}

\def\ot{\!\otimes\!}

\def\pp#1{\mathbb P^{#1}}
\def\ppp{\mathbb P}
\def\qq#1#2#3{q^{#1}_{{#2} {#3}}}

\def\ra{\rightarrow}
\def\rr#1#2#3#4{r^{#1}_{{#2} {#3}{#4}}}

\def\tbase{{\rm Base}\,}

\def\tdim{{\rm dim}\,}

\def\tker{{\rm ker}\,}

\def\trank{{\rm rank}\,}

\def\ttrace{{\rm tr}}

\def\up#1{{}^{({#1})}}
\def\upperp{{}^{\perp}}

\def\ww{\wedge}
\def\xx#1{x^{#1}}

\begin{document}

\title{On the projective geometry of rational homogeneous varieties}
\author{ J.M. Landsberg and Laurent Manivel}
\date{April 2002}

\begin{abstract}We determine the varieties of linear spaces on rational
homogeneous varieties, provide explicit geometric models for these spaces,
and establish basic facts about the local differential geometry of
rational homogeneous varieties.
\end{abstract}

\maketitle

\newcommand{\norm}[1]{\lVert#1\rVert}

\theoremstyle{plain}
\newtheorem{theo}{Theorem}[section]
\newtheorem{coro}[theo]{Corollary}
\newtheorem{lemm}[theo]{Lemma}
\newtheorem{prop}[theo]{Proposition}

\theoremstyle{definition}
\newtheorem{rema}[theo]{Remark}
\newtheorem{exam}[theo]{Example}
\newtheorem{defi}[theo]{Definition}

\section{Introduction}

Let $G$ be a complex simple Lie group and $P_S$ a parabolic subgroup
corresponding to a subset $S$ of nodes of the Dynkin diagram (so that a maximal parabolic subgroup
is defined by a single root).  Then
$G/P_S$ has a minimal homogeneous embedding in the projective space of
the highest weight module of $G$ corresponding to the weight
$\lambda= {\sum_{i\in S}\omega_i}$, where $\omega_i$ is the $i$-th
fundamental weight. We study the
local differential geometry of 
the embedded variety
$G/P_S\subset \BP V_{\lambda}$ and  the projective linear subspaces on
$G/P_S\subset \BP V_{\lambda}$. 

 We
describe the varieties parametrizing  such linear
 spaces in \S 4- 6. In most cases (those of \lq\lq non-short roots\rq\rq )
    the parameter 
varieties are determined in terms of Dynkin diagram data 
as explained in \S 4. (See in particular Theorem \ref{kplanes1}.)
The exceptional (exposed short root) cases 
are determined by use of explicit models in \S 5 for the
case of classical groups and \S 6 for the exceptional groups.
In all cases, each connected component of the variety
of linear spaces on a $G/P$
is quasi-homogeneous; more precisely, it is the union of  a finite number 
of $G$-orbits.

The case of unirulings by lines was studied in
\cite{coco} by means of Tits buildings. 
Our approach is by means
of projective differential geometry. 
This method is well suited because the
variety $G/P_S$ is homogeneous and
in particular cut out by quadrics, so the varieties of linear spaces on it
are determined by second order data at a point   $x\in G/P_S$.

In \S 2 we establish   basic connections between local differential
geometry and representation theory.  We study the semisimple
 part $H$ of $P_S$, which   fixes the point  
 $x=[v_{\lambda}]$
corresponding to the highest weight line  in its action on the
tangent space $T_xG/P_S$. As an $H$-module, $T_xG/P_S$ decomposes
into a direct sum of {\it generalized minuscule} $H$-modules.
If $S=\{\alpha\}$ where $\a$ is  a non-short root, the
space of tangent directions to lines is a {\it minuscule variety}
of $H$ (an irreducible, minimally homogeneously embedded Hermitian symmetric
space of $H$).

We study minuscule varieties in \S 3 and prove our main result
on their infinitesimal geometry in \S 4.

\smallskip

This is the first paper in  a series 
\cite{LM1,LM2, LMtrial, LMseries, LMpop}
 establishing new relations between 
the representation theory of complex simple Lie groups and the
algebraic and differential geometry of their homogeneous varieties.
The surprising connection between secant varieties and prolongations
developed in this paper 
is exploited in the sequels.
\footnote{ {\em Keywords}: flag variety, fundamental forms, 
linear spaces, octonion. {\em AMS classification}: 14M15, 20G05.}

\medskip\noindent {\em Acknowledgements}: We thank J. Wolf,
J.-M. Hwang and D. Snow for useful conversations, and an anonymous referee
for the simplified proof of Theorem 2.3.

\section{Under the microscope}
In this section we establish the basic connections between
differential invariants of homogeneous varieties and
representation theoretic data.
In \S\ref{ff} we  review the   projective fundamental forms of an
arbitrary projective variety $X^n\subset\pp\na$ and establish a connection between
secant varieties and fundamental forms.
In \S\ref{osc} we express the fundamental forms of homogeneously embedded
homogeneous varieties $X=G/P\subset \BP V_{\lambda}$ in terms of 
the universal envelopping algebra $U(\fg)$.
In \S\ref{decomp} we discuss the $P$-module
structure on $T_xX$, introduce an important class of homogeneous
varieties, the {\it generalized minuscule varieties} and explain
their role in the study of   fundamental forms of  rational homogeneous
varieties.

\subsection{Fundamental forms of projective varieties}\label{ff}

\subsubsection{Notation}
We let $V=\BC^{m}$ and $\BP V$  the corresponding
projective space. If $Y\subset \BP V$ is a set, we let $\hat Y\subset V$
denote the corresponding cone in $V$.
If $v\in V$, we let $[v]\in \BP V$ denote the
corresponding point in projective space.  

Let $X^n\subset \BP V=\pp\na$ be a projective variety of dimension $n$, and let
$x\in X$ be a smooth point. We let $T_xX$ denote the (intrinsic)
Zariski tangent space to $X$ at $x$, $\tilde T_xX\subset \pp\na$
denote the embedded tangent projective space (the $\pp n\subset\pp\na$
that best approximates $X$ at $x$), and $\hat T_xX\subset V =
\widehat{\tilde T_xX}$. We have the relation
$T_xX=\hat x^*\ot (\hat T_xX)/\hat x$
  and we  also  have, for
any $p\in V$ with $x=[p]$, $\hat T_xX=T^{affine}_p\hat X$, the affine tangent space at $p$.

We let $N_xX=T_x\BP V/T_xX$ denote the normal space of $X$ at $x$.

We use the ordering of the roots as in \cite{bou}.

\subsubsection{Fundamental forms in coordinates}
Let $X^n\subset\pp\na$ be a projective variety, and let
$x\in X$ be a smooth point. Take local linear coordinates
$(\xx 1\hd\xx n, \xx {n+1}\hd\xx\na )$
adapted to $x$, which means that they are centered at $x$ and
that $x^{n+1}=\cdots =x^{n+a}=0$ are equations of the embedded tangent space of $X$ at $x$. 
Write $X$, locally in the complex topology, as a graph
$$
\xx\mu = \sum_{1\leq \a,\b\leq n}\qq\mu\a\b\xx\a\xx\b + \sum_{1\leq \a,\b,\gamma \leq n}
\rr\mu\a\b\g\xx\a\xx\b\xx\g +\hdots ,
$$
where $n+1\leq\mu\leq n+a$.
The geometric information in the series (that is, information
independent of choice of adapted coordinates) can be encoded
in a series of tensors, the simplest of
which is the {\it projective second fundamental form}
$$\fff 2{X,x}=\sum_{\substack{1\le\a,\b\le n, \\ n+1\le \mu\le n+a}}
\qq\mu\a\b d\xx\a\circ d\xx\b \ot\frac{\partial}{\partial
\xx\mu}\in S^2T^*_xX\ot N_xX.$$ If $x$ is a general point, $\fff 2{X,x}$
even contains information about the global geometry of $X$,
see \cite{gh}, \cite{lanrr}. It is useful to consider the second fundamental 
form as a system of quadrics $\ii :=\ppp (\fff 2{X,x}(N^*_xX))\subseteq \ppp S^2T^*_xX$ 
parametrized by $N^*_xX$, and $\tbase\ii\subset\ppp T_xX$, their common zero locus.

We let $\hat T\up 2_xX=\hat T_xX +\cO (1)_x\ot\fff 2{X,x}(N^*_xX)\subset V$,
the {\it second osculating space to $X$ at $x$},
and $N_2=N_{2,x}X=\cO(-1)\ot (\hat T_x\up 2X/\hat T_xX)$.

More generally,   the {\it $k$-th projective fundamental}
form of $X$ at $x$ is  a tensor
$$
\fff k{X,x}\in S^kT^*_xX\ot N_{k,x}X
$$
where $N_k=N_{k,x}X=\cO(-1)_x\ot (\hat {T_x}\up kX/\hat {T_x}\up{k-1}X)$
and $\hat {T_x}\up kX=\hat {T_x}\up{k-1}X+ \cO(-1)_x\ot\ff{k }(S^{k }T_xX)$
is the {\it $k$-th osculating space to $X$ at $x$}.
To define $\fff k{X,x}$  
one can use the same definitions as one does
for the Euclidean fundamental forms, either in coordinates
or as the derivatives of successive Gauss mappings (see \cite{lanrr}).
Note that the osculating spaces
determine  a flag of $V$,
$$0\subset \hat x\subset \hat T_xX \subset
 \hat T_x\up 2X \subset ... \subset\hat T_x\up f =   V. $$

\smallskip
More generally, given a mapping $\phi : Y\ra \ppp V$,
one   defines its fundamental forms $\fff k{\phi}$ in the same
manner. $\fff 2{\phi,x}$ quotiented by  
  $\tker \phi_{*x}$ is isomorphic to 
the second fundamental form of the image, $\fff 2{\phi (Y),\phi (x)}$.
See    \cite{lansec} for details.

In what follows, we   slightly abuse notation by
ignoring twists
by the line bundles $\cO(j)$, which  will
 not matter as  we  study fundamental forms only at some
fixed base point. 
We let $|\fff k{X,x}|\subset \ppp S^kT^*_xX$ denote  
$\PP(\ff k_{X,x}(N^*_{k,x}X))$  and $\tbase |\fff k{X,x}|\subset\ppp T_xX$ denote its 
  base locus.
 
\smallskip

\subsubsection{Prolongation}
Let $V $ be  a vector space,   let $A\subset  S^dV^*$ be a linear 
subspace, and let
$$A\up l : =    (A\ot S^lV^*) \cap   S^{d+l}V^*,$$
 the {\em $l$-th  prolongation}
of $A$. Here the inclusion $S^{d+l}V^*\hookrightarrow S^dV^*\otimes S^{l}V^*$ is dual to the 
multiplication map $S^dV\otimes S^lV\lra S^{d+l}V$. If $P\in S^{d+l}V^*$, we still denote
by $P(v_1,\ldots, v_{d+l})$ its polarization. Then $P\in A\up l$ if and only if for all $v\in V$, 
the degree $d$ polynomial $w\mapsto P(v,\ldots ,v,w,\ldots ,w)$ belongs to $A$.  
The notation is such that $A\up 0=A$. 

Let  $Jac(A):=   \{ v\lrcorner P | v\in V, P\in A \}\subseteq
  S^{d-1}V^*$,  the {\em Jacobian space of
$A$}. Then  $A\up 1 =    \{ P\in  S^{d+1}V^* | Jac(P)\subset A\}$.

\smallskip

A basic fact about fundamental forms,
due to Cartan (\cite{car}, p. 377) (and rediscovered in \cite{gh}), is that 
if $x\in X$ is a general point, then the {\em prolongation
property} holds at $x$:
$$ |\fff k{X,x}|\subseteq |\fff {k-1}{X,x}|\up 1.$$ 

 A  geometric consequence is as follows. 
Define the {\it $k$-th secant variety} $\s_k(Y)$ of a projective
variety $Y\subset\PP^N$ to be the closure of the union of the linear spaces 
spanned by $k$ points of $Y$. The notation is such that $\s_1(Y)=Y$.

\begin{prop}\label{sec1}
Let $X^n\subset\pp\na$ be a  variety and $x\in X$ a general point. Then
for $k\ge 2$, 
$$\tbase |\fff k{X,x} |\supseteq \sigma_{k-1}(\tbase |\fff 2{X,x}|). $$
 \end{prop}

 Proposition  \ref{sec1}  is a consequence of the following lemma:
 
 \begin{lemm} Let $A\subset S^2V^*$ be a system of quadrics with base locus $\tbase (A)\subset \PP V$. 
Then $$\tbase (A\up{k-1})\supseteq \s_k(\tbase (A)).$$
Moreover, if $\tbase (A)$ is linearly non-degenerate, then for $k\geq 2$, 
$I_k(\s_k(\tbase (A))=0$, and if $A=I_2(\tbase (A))$, then 
$I_{k+1}(\s_k(\tbase (A))=A\up k$, where $I_d(Z)\subset S^dV^*$ is the 
component of the ideal of $Z\subset\ppp V$ in degree $d$.
\end{lemm}

\proof We prove the lemma for $k=2$, the generalization being clear. 
We first need to prove that 
any polynomial $P\in A \up 1 $ vanishes on $v=sx+ty$ for all 
$s,t\in \CC$ and $x,y\in B$, the cone over $\tbase (A)$. 
Since $P(x,x,.)=P(y,y,.)=0$, we have 
$$P(v) = P(v,v,v)=
s^3P(x,x,x) + 3s^2tP(x,x,y) + 3st^2P(x,y,y) + t^3P(y,y,y)=0.$$

Now, say $Q\in I_2(\s_2(\tbase(A))$. Then for
all $x,y\in B$ and $s,t\in\CC$, $Q(sx+ty)=0$,  which implies
$Q(x,y)=0$, which implies $Q=0$ since $\tbase (A)$ is non-degenerate. 

Finally, consider a polynomial $P\in  I_3(\s_2(B(A))$. Since $P$ 
vanishes on $v=sx+ty$ for all $x,y\in B$ and all $s,t\in\CC$, 
we have $P(x,x,y)=0$ for all $x\in B$, and all $y\in  B$, 
hence all $y\in V$ since $\tbase (A)$ is non-degenerate. 
Thus for all $y\in V$, 
$P(.,.,y)$ is a quadric vanishing on $\tbase (A)$, hence belongs to 
$A=I_2(\tbase (A))$. This means that $P$ is in $A\up 1$. 
 \qed 

\medskip
An elementary fact about projective varieties is that if 
$X^n\subset\pp\na$ is  a  variety whose ideal is generated 
in degree $\leq d$, and $L$ a linear space osculating to order
$d$  at a smooth point $x\in X$, then  $L\subset X$.
The ideal of  a  projective 
homogeneous variety is generated in degree two (see e.g.  
\cite{lich}), so if $X\subset\ppp V$ is homogeneous,  then $\tbase\ii$
is the set of tangent directions to lines on $X$ through $x$.
If $y\in \tilde T_xX\cap X$ then the line $\pp 1_{xy}$ is contained in $X$. 

\subsection{Osculating spaces of homogeneous varieties}\label{osc} 

Let $G$ be a simply connected complex semi-simple Lie group, $\fg$ its 
Lie algebra,   $\fg^\otimes$ the tensor algebra of $\fg$ and  $U (\fg )
=\fg^{\ot }/
\{x\ot y - y\ot x - [x,y]\mid x,y\in\fg\}$,
the universal envelopping algebra. $U(\fg)$
inherits a filtration from the natural   grading of $\fg^\otimes$, and the
associated graded algebra is the symmetric algebra of $\fg$. 
Fix a maximal torus $T$ and a Borel subgroup $B$ of $G$ containing
$T$. We adopt the convention that $B$ is generated by the positive
roots, and we write the corresponding root space decomposition of $\fg$ as
$$
\fg  =     \ft \oplus \bigoplus_{\alpha\in \Delta_+}(\fg_{\alpha} \oplus
\fg_{-\alpha} ),
$$ 
where $\Delta_+$ denotes the set of positive roots.  

Let $V_{\l}$ be an irreducible $G$-module with highest weight $\l$, and
$v_{\l}\in V_{\l}$ a highest weight vector. 
The induced action of $\fg$ extends to the universal envelopping algebra,
inducing a  filtration of $V_{\l}$ whose $k$-th term is  
$$
V_{\l}^{(k)} =      
U_k(\fg )v_{\l}.
$$

Let $x =    [v_{\l}] \in \BP V_{\l}$  and 
let $X  =    G/P\subset \PP V_{\l}$ be its $G$-orbit.
 Here $P$ 
is the stabilizer of $x $, it is a parabolic subgroup of $G$. The
tangent bundle $TX$ is a homogeneous bundle and we identify $T_xX$ with
the associated
 $P$-module $\fg /\fp$. The osculating spaces and the fundamental forms
of $X$ have a simple representation-theoretic interpretation: 

\begin{prop}
Let $X=   G/P \subset\ppp V_{\l} $
be a homogeneous variety  with base point $x = [v_{\l}]$.
Let $\hat T^{(k)}_{x }X$ denote the cone over the $k$-th osculating
space at
$x $ and let $N_k=    \hat T^{(k)}_{x }X/\hat T^{(k-1)}_{x }X$ be the
 $k$-th
normal space twisted by ${\mathcal O}(-1)$.    Then
$$
\hat T^{(k)}_{x }X =    V_{\l}^{(k)}, \qquad
N_k =    V_{\l}^{(k)}/V_{\l}^{(k-1)}.$$
  Moreover, there is a commutative diagram   
$$ \begin{array}{rcc}
S^k\fg &  =     & U_k(\fg )/U_{k-1}(\fg )\\  
\downarrow\hspace*{3mm} & & \downarrow \\
  \FF^k_{X,x } : S^kT_xX & \longrightarrow & N_k, 
\end{array}
$$
where the bottom horizontal map is  
the   $k$-th  fundamental form at $x$. 
\end{prop}

\proof The diagram above is   the $k$-th fundamental form of the
mapping $\phi : G\ra  \ppp V$ at $e\in G$,
where $\phi (G)=X$. 
 \qed

\medskip  $V_{\l}^{(k)}$ has a natural $P$-module structure.
 Thus the osculating spaces of $X$ at $x$ correspond to the increasing
filtration of
$P$-modules  
$$
0\subset \hat x \subset \hat
T_{x }X  =     V_{\l}^{(1)} \subset V_{\l}^{(2)} \subset \cdots \subset
V_{\l}^{(f)}  =      V_{\l}.
$$
Our next goal is to understand the first quotient of this  filtration,
namely the structure of $T_xX$ as a $P$-module.

\subsection{Decomposing the tangent space} \label{decomp} 

Let $P=P_{\a_i}$ denote the maximal parabolic subgroup of $G$ 
corresponding to the simple root $\a_i$. Let $P = LP^u$ be a Levi 
decomposition of $P$, where $P^u$ is unipotent, $L$ is reductive and contains 
the maximal torus $T$. If $\a$ is a positive root, let $\a  =    
\sum_jm_j(\a)\a_j$ be its decomposition in terms of simple roots. Let 
$\Delta_X  = \{\alpha\in \Delta_+ | m_i(\a)>0\}$. We have the root space 
decompositions 
$$ \begin{array}{lcl}
\fp &  =     & \ft \oplus (\bigoplus_{\alpha\in \Delta_+}\fg_{\alpha} ) 
\oplus (\bigoplus_{\alpha\in \Delta_+\backslash \Delta_X}\fg_{-\alpha}),\\
\fl &  =     &
\ft \oplus \bigoplus_{\alpha\in \Delta_+\backslash\Delta_X} (\fg_{\alpha}
 \oplus \fg_{-\alpha}), \\
\fp^u &  =    & \bigoplus_{\alpha\in \Delta_X}\fg_{\alpha}. \end{array}$$

\begin{prop}\label{tangent} 
Let $G$ be simple, let $\tilde\a$ be the highest root 
of $\fg$, let $\a_i$ be a simple positive
root, and let $P  =     P_{\a_i}$ be the associated 
maximal parabolic subgroup. For $1\le k\le m_i(\tilde\a)$, let
$$\fs_k  =    \bigoplus_{m_i(\a)\ge -k}\fg_{\alpha}.  $$
This defines an increasing filtration of $\fg /\fp$ by $P$-submodules.
The quotients 
$$T_k  =    \bigoplus_{m_i(\a)=k}\fg_{-\alpha}  $$
are irreducible $P$-modules.\end{prop} 

\proof The fact that each $\fs_k/\fs_{k-1}$ is a $P$-module is clear. 
The irreducibility of  $T_k$ is a special case of  \cite{Wolf}, 8.13.3
(which is attributed to Kostant). \qed 
 
\medskip
The irreducibility of $T_k$ implies  that the set
$\{ \a\in \Delta_+ \mid m_i(\a )=   k\}$  has a unique minimal element
which we denote by $-\phi_k$ when we consider the root as a weight of $T_k$. 
In particular, the highest weight of $T_1$   is  
$$ \phi_1  =     -\a_i  =     -\sum_jn(\a_i,\a_j)\omega_j,  $$
where $n(\a_i,\a_j)$ denotes the entries of the Cartan matrix.  This
weight is easy to read directly on the Dynkin diagram of $G$. 
Let $H$ denote the semi-simple part of $L$. As an $H$-module, the
filtration of $T_xX$ into irreducible $P$-modules becomes a
 direct sum decomposition into irreducible $H$-modules.
Note that $L$ has a one dimensional center  and that the Lie algebra of
$H$ is 
$$
\fh  =     \ker\o_i \oplus 
\bigoplus_{\alpha\in \Delta_+\backslash\Delta_X}
 (\fg_{\alpha} \oplus \fg_{-\alpha}),
$$
where $\ker\o_i\subset\ft$.
The Dynkin diagram of $H$ is therefore deduced from that
 of $G$ by   suppressing the node corresponding to the simple root
$\a_i$.  In particular, we conclude:

\begin{prop} Let $X=   G/P$ be a homogeneous variety
with $P$ a maximal parabolic and let $H$ be the
semi-simple part of $P$. Then $T_1$, the first irreducible
component of $T_xX$ as an $H$-module,  is  
  obtained by marking the nodes of ${\mathcal D}(H)$ adjacent to
the node from ${\mathcal D}(G)$ that was removed. A node $\b$ is
given   multiplicity two (resp. three) if there is an arrow emanating
from $\a$ towards $\b$ with a double (resp. triple) bond.
\end{prop}

The above observations can be found in  \cite{dynkin}.  
 
\begin{defi} A fundamental weight $\o_i$ is {\em minuscule} if 
the Weyl group acts transitively on the set of weights of the 
corresponding fundamental representation. 

In an irreducible root
system, a fundamental weight $\o_i$ is {\em cominuscule} 
if the highest root has coefficient one on $\a_i$. In a reducible
root system, a weight is cominuscule if it is a sum of cominuscule 
fundamental weights, one for each irreducible factor of the root
system. \end{defi}

The relation between these two notions is as follows. In the
irreducible case, $\o_i$
is minuscule if and only if in the dual root system, the highest
root has coefficient one on the coroot $\check\a_i$ (\cite{bou},
Chap. 8). 

Geometrically, when $G$ is simple, the weight $\o_i$ is cominuscule
exactly when $G/P_i$ admits the structure
of an irreducible Hermitian  symmetric space whose
automorphism group is locally isomorphic to $G$
(we call $G/P_i$ a $G$-Hermitian symmetric space). This was 
pointed out by Kostant in \cite{kos}. We   use the following 
definition (be careful that minuscule varieties are in correspondence
with cominuscule weights, not minuscule weights!): 

\begin{defi} A {\em $G$-minuscule variety} $X=G/P\subset \PP V$ is a 
$G$-Hermitian (not necessarily irreducible) symmetric space in its 
minimal homogeneous embedding. A {\it generalized minuscule variety} is 
a Hermitian symmetric  space $X=G/P\subset\ppp V$ in some $G$-homogeneous embedding, 
but the automorphism group of $X$ need not be locally isomorphic to
$G$, and the embedding need not be minimal. \end{defi}

\begin{prop}\label{y1.2}
Let  $X =  G/P$ with $G$  simple, and $P=P_{\a_i}$ a maximal
subgroup with semi-simple part $H$. If
$\a_i$ is not short,
then the closed $H$-orbit $Y_1\subset \PP T_1$ is an 
$H$-minuscule variety. 
\end{prop}

See \ref{moreg} below for a more general statement.

\proof 
The observations above imply that
$Y_1=H_1/Q_1\times\cdots\times H_r/Q_r$, where $H_1,\ldots, H_r$
are simple Lie groups whose Dynkin diagrams are the branches from
$\a_i$ of ${\mathcal D}(G)$, and $Q_1,\ldots ,Q_r$ are maximal
parabolic subgroups defined for each branch by the node adjacent 
to $\a_i$. Since the root system of each $H_k$ is formed by 
the roots of $G$ with support on the corresponding branch, we just
need to prove that if a root $\b$ of $G$ has coefficient $m_i(\b)=0$
on $\a_i$, its coefficient on an adjacent root $\a_k$ cannot exceed
one. But this follows immediately from the equality $n(\b,\a_i)=
\sum_jm_j(\b)n(\a_j,\a_i)$ (where the integers $n(\a_j,\a_i)$ are
non-positive, and negative exactly when $\a_j$ is adjacent to $\a_i$),
and the fact that, since $\a_i$ is not short, $|n(\b,\a_i)|\le 1$. 
The minimality of the embedding of $Y_1$ in $\PP T_1$ similarly 
follows from
the fact that $n(\a_j,\a_i)=-1$ for $\a_j$ connected to $\a_i$. \qed

\medskip We can say slightly more when $G$ is simply laced. 

\begin{prop} Let $G$ be a simple Lie group
of type $A,D$ or $E$,
 let $P$ be a maximal parabolic subgroup, let $T_1\subset T_xG/P$
be as in \ref{tangent}. Then $T_1$ is a minuscule $H$-module. \end{prop}

\proof The weights of $T_1$ are, by definition, the opposites of 
the roots $\b$ such that $m_i(\b)=1$. If $G$ is simply laced, 
these roots all have the same length and lemma \ref{slem},  
which was communicated to us
anonomously, shows that the Weyl group of $H$ acts 
transitively on them. \qed

\medskip The above discussion   can easily be extended to homogeneous spaces $X=G/P$
with $P$ not necessarily maximal. Suppose that  $P=P_S$ is the 
parabolic subgroup generated  by the complement of a set $S$ of simple 
positive roots.  Then  there is an irreducible component of the
$L$-module $T=T_x(G/P)$ for each   choice of the coefficients of the positive 
roots on these simple roots. If we choose such a family of
coefficients $a = (a_i)_{i\in S}$,  and let
$$\fs _a =    \bigoplus_{\a\in\Delta_+,\;m_i(\a)\ge -a_i}\fg_{\a}
\quad {\rm and} \quad
T_a =    \bigoplus_{\a\in\Delta_+,\;m_i(\a) =  a_i}\fg_{-\a},$$ 
then $\fs _a$ is a $P$-submodule of $\fg/\fp$, and $T_a$ is an
irreducible $L$-submodule of $\fs_a$. An important difference 
with the case of maximal parabolics is that the incidence relations 
between the non-zero $\fs_a$'s is no longer a simple chain of inclusion, 
but defines a partial order.  

\smallskip
Let $\e_i$ be the family of coefficients $a_i=1$, and $a_j=0$ for $j\in S-i$.
The analogues of $T_1$ in the maximal case are the $H$-modules
$T_{\e_i}$ (note that only the irreducible factors of $H$ corresponding 
to the branches of $\cD\backslash S$ connected to $\a_i$ act non 
trivially on $T_{\e_i}$; we denote their product by $H_{\e_i}$).
Again we need to know whether $\a_i$ is short or not, but this condition
is relevant only with respect to a subdiagram of ${\mathcal D}(G)$. 

\begin{defi} We call $\a\in S$ an {\it exposed short root} if the
connected component of $\a$ in $\cD\backslash (S\backslash \a)$
contains a root longer than $\a$, i.e., if an arrow in 
$\cD\backslash (S\backslash \a)$ points towards $\a$. \end{defi}

\begin{prop}\label{moreg} Notations    as above. 
Let $Y_1^i\subset \PP T_{\e_i}$ be the closed orbit. 
Then $Y_1^i$ is  a generalized minuscule variety.
Moreover, it is a $H_{\e_i}$-minuscule variety, 
except for situations equivalent to the following cases: 

1. $ C_n/P_k$ for $k<n$. Here $Y_1= Seg(\PP^{k-1}\times\PP^{2n-2k-1})$, 
$H= SL_k\times Sp_{2n-2k}  \subsetneq SL_k\times SL_{2n-2k}$.

2. $ C_n/P_n$. Here $Y_1=v_2(\pp{n-1})$, which is $A_{n-1}$-minuscule
but not in its minimal embedding.

3. $F_4/P_4$. Here $Y_1= B_3/P_{3}$, a six-dimensional quadric.

4. $G_2/P_2$. Here $Y_1= v_3(\pp 1)$ is the twisted cubic, which is
$A_2$-minuscule, but not in its minimal embedding. 
\end{prop}

Let $X=G/P$ with $P$ maximal, let $H$ the semi-simple
part of $P$.  We obtain a splitting $T_xX=\oplus_pT_p$, with
each $T_p$ an irreducible $H$-module. Let $Y_p\subset \BP T_p$
denote the closed orbit.

 \begin{prop} \label{ratcurv}    The closed orbit 
$Y_{p}$ is contained in $\tbase |\fff {p+1}{X,x} |,$ 
and there is a rational normal curve in $X$ of degree at most $p+1$, 
passing through $x$ with tangent vector in $Y_{p} $.
\end{prop}

This proposition indicates that it is possible to study the 
$G$-homogeneous rational curves on $G/P$ of  degree
greater than one using the methods we use to study lines on $G/P$.

\proof 
Let $P=P_{\a_j}$ and let $\b$ be such that
$m_j(\b)=p$.
Let $X_{\b}\in \fg_{-\b}$.
Let $v\in V$ be a highest weight vector and $x=[v]$.
Then $X_{\b}v\in \hat Y_{p}\subset T_p$.
By \cite{dix}, lemme 7.2.5, $X_{\b}^{p+1}v =0$ so the 
rational curve 
$  \exp (tX_{\b})v$
 is
contained in $X$
and   is of degree at most $p$.\qed

\medskip
These propositions stress the importance of minuscule varieties 
in our study. The next section is devoted to their properties.

\section{Minuscule varieties}

We explicitly describe the tangent and normal spaces to 
minuscule varieties  $X=G/P$
in \S\ref{minutan} and \S\ref{minunormal}
as $H$-modules, where $H$ is the semi-simple part of the Levi factor
of $P$. In \S\ref{minuthm} we state an prove our main theorem that 
determines
the fundamental forms of minuscule varieties. In \S\ref{algstrs}
we remark on some interesting complexes obtained from the
normal spaces of minuscule varieties.

\subsection{Their tangent spaces}\label{minutan}

We summarize characterizations and tangent space
structures of minuscule varieties:
Let  $G$ be a  simple Lie group and $P  = P_{\a_i}$ a maximal parabolic
subgroup. Let $G/P\subset\ppp V$ be the minimal homogeneous embedding. 
The following are equivalent:

1.  $m_i(\tilde\a) =    1$ (the highest root $\tilde\a$  has coefficient
one on the simple root $\a_i$),

2.  $\fp^u$ is an abelian subalgebra of $\fg$,

3. $T=T_{[e]}(G/P)$ contains no $P$-invariant submodule,

4.  $G/P$ admits an irreducible Hermitian symmetric metric
with local holonomy $G$ induced from a Fubini-Study metric on $\ppp V$,
and the embedding to $\ppp V$ is the smallest such embedding.

\smallskip
Here is a table of the $G$-minuscule varieties: there 
are four infinite series and two exceptional spaces. 

\medskip
\begin{center}\begin{tabular}{||c||c|c|c|c||} \hline  
Name & ${\scriptstyle Grassmannian}$ &
${\scriptstyle Quadric}$ & ${\scriptstyle Lagrangian\, Grassm.}$ & 
${\scriptstyle Quadric}$   \\  \hline  
Notation & $G(k,n+1)$ & $\QQ^{2n-1}$ & $G_{Lag}(n,2n)$ & 
$\QQ^{2n-2}$ \\ 
$G$ & $A_n$ & $B_n$  & $C_n$ & $D_n$ \\  
$\omega$ & $\o_k$ & $\o_1$ & $\o_n$ & $\o_1$ \\ 
${\mathcal D}(G)$ & \setlength{\unitlength}{2.5mm} 
\begin{picture}(10,3)(0.5,-1)
\multiput(0,0)(2,0){6}{$\circ$}
\multiput(0.55,0.35)(2,0){5}{\line(1,0){1.55}}
\put(6,0){$\bullet$}\end{picture} & 
\setlength{\unitlength}{2.5mm} \begin{picture}(10,3)(0.5,-1)
\multiput(0,0)(2,0){6}{$\circ$}
\multiput(0.55,.35)(2,0){4}{\line(1,0){1.55}}
\multiput(8.55,.25)(0,.2){2}{\line(1,0){1.55}} \put(8.75,0.02){$>$}
\put(0,0){$\bullet$}\end{picture} & 
\setlength{\unitlength}{2.5mm} \begin{picture}(10,3)(0.5,-1)
\multiput(0,0)(2,0){6}{$\circ$}
\multiput(0.55,.35)(2,0){4}{\line(1,0){1.55}}
\multiput(8.55,.25)(0,.2){2}{\line(1,0){1.55}} \put(8.85,0.03){$<$}
\put(10,0){$\bullet$}\end{picture} &
\setlength{\unitlength}{2.5mm} \begin{picture}(10,3)(0.5,-1)
\multiput(0,0)(2,0){5}{$\circ$}
\multiput(0.55,.4)(2,0){4}{\line(1,0){1.55}} 
\put(8.5,.35){\line(2,-1){1.55}}
\put(8.5,.35){\line(2,1){1.55}}
\put(10,.9){$\circ$}
\put(10,-.9){$\circ$}
\put(0,0){$\bullet$}\end{picture}  \\ 
$H$ & $A_{k-1}\times A_{n-k}$ & $B_{n-1}$  & $A_{n-1}$  & 
$D_{n-1}$   \\   
$\phi_1$ & $\o_{k-1}+\o_{k+1}$ & $\o_1$ & $2\o_{n-1}$ & $\o_1$ 
\\ 
${\mathcal D}(H)$ & \setlength{\unitlength}{2.5mm}
\begin{picture}(10,3)(0.5,-1)
\multiput(1,0)(2,0){2}{$\circ$}
\multiput(9,0)(2,0){1}{$\circ$}
\multiput(1.55,0.35)(2,0){2}{\line(1,0){1.55}}
\multiput(7.55,0.35)(2,0){1}{\line(1,0){1.55}}
\put(5,0){$\bullet$}\put(7,0){$\bullet$} \end{picture} &
\setlength{\unitlength}{2.5mm}
\begin{picture}(10,3)(0.5,-1)
\multiput(2,0)(2,0){5}{$\circ$}
\multiput(2.55,.35)(2,0){3}{\line(1,0){1.55}}
\multiput(8.55,.25)(0,.2){2}{\line(1,0){1.55}} \put(8.8,-.05){$>$}
\put(2,0){$\bullet$}\end{picture} & 
\setlength{\unitlength}{2.5mm}
\begin{picture}(10,3)(0,-1)
\multiput(0,0)(2,0){5}{$\circ$}
\multiput(0.55,.35)(2,0){4}{\line(1,0){1.55}}
\put(8,0){$\bullet$}\put(8,.85){${\scriptstyle 2}$} \end{picture} & 
\setlength{\unitlength}{2.5mm} \begin{picture}(10,3)(1,-1)
\multiput(2,0)(2,0){4}{$\circ$}
\multiput(2.55,.35)(2,0){3}{\line(1,0){1.55}} 
\put(8.5,.35){\line(2,-1){1.55}}
\put(8.5,.35){\line(2,1){1.55}}
\put(10,.9){$\circ$}
\put(10,-.9){$\circ$}
\put(2,0){$\bullet$}\end{picture} \\   
$T$ & $E^*\ot Q$ & $E^*\ot (E^{\perp}/E)$ &  $S^2Q$ & $E^*\ot
(E^{\perp}/E)$
\\ \hline \end{tabular}\end{center}\medskip

\begin{center}\begin{tabular}{||c||c|c|c||} \hline  
Name & ${\scriptstyle Spinor\, variety}$ & 
${\scriptstyle Cayley\, plane}$ & ?? \\  \hline
Notation &  ${\mathbb S}_n$ & ${\mathbb O }{\mathbb P }^2$ &
$G_{\o}(\BO^3, \BO^6)$ \\
$G$ & $D_n$ & $E_6$  & $E_7$ \\ 
$\o$ & $\o_n$ & $\o_1$ & $\o_7$ \\ 
${\mathcal D}(G)$  &
\setlength{\unitlength}{2.5mm} \begin{picture}(10,3)(0.5,-1)
\multiput(0,0)(2,0){5}{$\circ$}
\multiput(0.55,.35)(2,0){4}{\line(1,0){1.6}} 
\put(8.55,.35){\line(2,-1){1.55}}
\put(8.55,.35){\line(2,1){1.55}}
\put(10,.9){$\bullet$}
\put(10,-1.1){$\circ$}\end{picture} &
\setlength{\unitlength}{2.5mm}
\begin{picture}(9,3)(0,-1.8)
\multiput(0,0)(2,0){5}{$\circ$}
\multiput(0.55,.35)(2,0){4}{\line(1,0){1.6}} \put(0,0){$\bullet$}
\put(4,-2){$\circ$}
\put(4.35,-1.4){\line(0,1){1.5}}\end{picture} & 
\setlength{\unitlength}{2.5mm}
\begin{picture}(11,3)(0.8,-1.8)
\multiput(0,0)(2,0){6}{$\circ$}
\multiput(0.5,.3)(2,0){5}{\line(1,0){1.6}} \put(10,0){$\bullet$}
\put(4,-2){$\circ$}
\put(4.35,-1.4){\line(0,1){1.5}}\end{picture} \\ 
$H$ & $A_{n-1}$  & $D_5$  & $E_6$  \\ 
$\phi_1$ & $\o_{n-2}$ & $\o_4$ & $\o_6$ \\ 
${\mathcal D}(H)$ & \setlength{\unitlength}{2.5mm} 
\begin{picture}(10,3)(0.5,-1)
\multiput(0,0)(2,0){5}{$\circ$}
\multiput(0.55,.35)(2,0){4}{\line(1,0){1.6}}
\put(8.55,.35){\line(2,-1){1.55}}
\put(8,0){$\bullet$}
\put(10,-1.1){$\circ$}\end{picture} &
\setlength{\unitlength}{2.5mm}
\begin{picture}(9,3)(1,-1.8)
\multiput(2,0)(2,0){4}{$\circ$}
\multiput(2.55,.35)(2,0){3}{\line(1,0){1.6}} \put(2,0){$\bullet$}
\put(4,-2){$\circ$}
\put(4.35,-1.4){\line(0,1){1.5}}\end{picture} & 
\setlength{\unitlength}{2.5mm}
\begin{picture}(9,3)(0.8,-1.8)
\multiput(0,0)(2,0){5}{$\circ$}
\multiput(0.55,.35)(2,0){4}{\line(1,0){1.6}} \put(8,0){$\bullet$}
\put(4,-2){$\circ$}
\put(4.35,-1.4){\line(0,1){1.5}}\end{picture} \\ 
$T$ & $\La 2 E^*$ & ${\mathcal S}^+$ &  $\JO$ \\ \hline
\end{tabular}\end{center}\medskip

Here $E$ and $Q$ are
the tautological and quotient vector bundles on  the Grassmannian
or their pullbacks to the varieties in question. 
${\mathcal S}^+$ is the half spin representation of $D_5$,
and $\JO$ is the space of $3\times 3$ $\BO$-Hermitian matrices,  
  the representation $V_{\lo 1}$ for $E_6$ (see \S 6.2 for
details).
$G_{\o}(\BO^3, \BO^6)$ may be interpreted as the space of $\BO^3$'s in $\BO^6$
that are null for an  $\BO$-Hermitian symplectic form, see \cite{LM1}.

\subsection{The strict prolongation property}\label{minuthm}

We prove our main theorem on the infinitesimal geometry of
minuscule varieties.
 
\begin{theo}
Let $X =    G/P_{\a_i}\subset \PP V_{\o_i}$ be a  minuscule variety
and $x\in X$.   Then for $k\ge 2$, 
$$
|\FF^{k+1}_{x,X}|=    |\FF^2_{x,X}|\up {k-1}.
$$
 \end{theo}

\begin{rema} This result says that the leading terms of the Taylor series in
local  coordinates adapted to the filtration by osculating spaces,
are determined by the quadratic terms in an elementary manner. 
In \cite{LM1}, we show moreover that there are no terms in the Taylor series
except for the leading terms. (Minuscule varieties are the unique
homogeneous varieties having this property.)\end{rema}

\proof
Let $v =    v_{\o_i}\in V_{\o_i}$ be the highest weight vector,
and let $T=T_{[v]}X$.
  We denote by $R_k\subset S^kT$ the space of relations of
degree $k$, that is, the space of homogeneous polynomials $P_k$ of degree
$k$ in the $X_{\a}$, with $\a\in\Delta_X$, such that   $P_k.v\in \hat
T^{(k-1)}_{[v]}X$, the $(k-1)$-st osculating space. We have the 
following  commutative diagram, where
 horizontal middle long sequence and the vertical short sequences 
 are exact:

$$\begin{array}{ccccccccc}
 & & 0 & & 0 & & 0 & & \\
 & & \downarrow & & \downarrow & & \downarrow & & \\
 \cdots & \lra & R_{k-1}\ot\Lambda^2T & \lra & R_k\ot T & \lra
& R_{k+1} & \lra & 0 \\
 & & \downarrow & & \downarrow & & \downarrow & & \\
\cdots & \lra & S^{k-1}T\ot\Lambda^2T & \lra & S^kT\ot T & \lra
& S^{k+1}T & \lra & 0 \\
& & \downarrow & & \downarrow & & \downarrow & & \\
\cdots & \lra & N_{k-1}\ot\Lambda^2T & \lra & N_k\ot T & \lra
& N_{k+1} & \lra & 0 \\
&  & \downarrow & & \downarrow & & \downarrow & & \\
&  & 0 & & 0 & & 0 & & 
\end{array}$$

\begin{lemm}   $N_{k+1}^* =    N_k^{*(1)}$
for all $k\geq 2$ if and only if the relations are
  generated in degree two, that is, the map $R_k\ot T\ra R_{k+1}$ 
is surjective for all $k\ge 2$. \end{lemm}

\proof We first note that $N_{k+1}^* =    N_k^{*(1)}$ holds if and only if
the  sequence 
$$N_{k+1}^*\lra N_k^*\ot T^* \lra N_{k-1}^*\ot\Lambda^2T^*$$
is exact at the middle term.  This is because,
by definition, $N_k^*\up 1=   (N_k^*\ot T^*)\cap S^{k+1}T^*$
and   $S^{k+1}T^*$ is
the kernel of the map $S^kT^*\ot T^*\ra S^{k-1}T^*\ot \La 2T^*$.

A diagram chase, using the above 
 partially exact diagram, shows that the exactness of the dual 
sequence
$ N_{k-1}\ot\Lambda^2T\lra  N_k\ot T \lra N_{k+1}$
 is equivalent to the surjectivity of the map $R_k\ot T\ra R_{k+1}$. \qed

\medskip Now  we analyze the space of relations. 
By  \cite{dix} (Lemme 7.2.5 p. 225), 
the relations all come from the identities 
$$X_{\a_i}^2v =    0 \quad {\rm and} \quad 
X_{\b}v =    0\;\; {\rm for}\; m_i(\b) =    0.$$

More precisely, if $P_k$ is a homogeneous relation of degree $k$, 
there   exists an 
identity of the following kind in     $U(\fn )$
(where $\fn$ is the subalgebra of $\fg$ generated by positive root
vectors):
$$P_k+Q_{<k}+\sum_{m_i(\b) =    0}R_{\b}X_{\b}+SX_{\a_i}^2 =    0,$$
where $Q_{<k}$ is a polynomial of degree less than $k$ in the 
$X_{\a}$, $\a\in\Delta_X$, and the $R_{\b}$ and $S$ are polynomials 
in the $X_{\g}$, $\g\in\Delta_+$.  

Now we fix an ordered basis of $\fn$, beginning first with the $X_{\b}$,
$\b\neq\a_i$,
such that $m_i(\b)>0$, then $X_{\a_i}$, 
and then  continuing with the $X_{\g}$ for which $m_i(\g) =    0$. By the 
Poincare-Birkhoff-Witt theorem ( \cite{dix}, Th\'eor\`eme 2.1.11
p. 69),  the monomials in the 
$X_{\g}$ compatible with this order  form a basis of   $U(\fn )$. 

We    say that a polynomial expression in the $X_{\g}$,
$\g\in\Delta_+$, is {\em well-ordered}  if each of its monomials is
compatible with our ordered basis. We may suppose that in the 
identity above, all the polynomials $P_k$, $Q_{<k}$, $R_{\b}$
and $S$ are well-ordered. We may even suppose that the products
$R_{\b}X_{\b}$ are well-ordered, as if they are not,
reordering them  
gives a sum of expressions of the same type,
 since the 
space generated by the $X_{\g}$ for which $m_i(\g) =    0$ is  
stable under the Lie bracket. However, 
and this is the crucial point, we cannot suppose {\it a priori} that the 
product $SX_{\a_i}^2$ is also well-ordered.

The conclusion of this analysis is   that all relations 
appear in the following way: we first chose a well-ordered monomial 
$X_{\b_1}\cdots X_{\b_m}$, with $m_i(\b_1) =    \cdots  =    m_i(\b_m) =   
0$; we reorder its product with $X_{\a_i}^2$, which gives 
an  expression of the form:
$$
X_{\b_1}\cdots X_{\b_m}X_{\a_i}^2 =    \sum_{\g\d}c_{\g\d}X_{\a_i+\g}
X_{\a_i+\d}+CX_{2\a_i+\b_1+\cdots +\b_m}+
\sum_{m_i(\eta) =    0}U_{\eta}X_{\eta},
$$
where $U_{\eta}$ is some polynomial in the $X_{\g}$
and $C$ is a constant. 
We then multiply on the left by a monomial in the $X_{\b}$ 
with $m_i(\b)>0$ and reorder if necessary, then we make linear 
combinations, and finally, we only keep the homogeneous terms
of maximal degree in the resulting expression. 

This doesn't seem very enlightening, but since $X$ is a minuscule variety,
if  $m_i(\b) =    m_i(\g)=    1$, then $X_{\b}$
and $X_{\g}$   commute. So   the above relation   simplifies 
to an expression of the form 
$$
X_{\b_1}\cdots X_{\b_m}X_{\a_i}^2 =    \sum_{\g\d}c_{\g\d}X_{\a_i+\g}
X_{\a_i+\d}+\sum_{m_i(\eta) =    0}U_{\eta}X_{\eta}.
$$
Moreover, the relations are then obtained by multiplying the 
sums $\sum_{\g\d}c_{\g\d}X_{\a_i+\g}X_{\a_i+\d}$ by monomials in the
 $X_{\a_i+\eta}$, which need no reordering; and finally, the 
resulting expression is necessarily homogeneous, since we can
assume that all its monomials have the same total weight.

This means in particular that all
relations are deduced from the degree two relations 
$$\sum_{\g\d}c_{\g\d}X_{\a_i+\g}X_{\a_i+\d} = 0$$
by simple polynomial multiplication in $T$. Thus the 
maps $R_2\ot S^{k-1}T\ra R_{k+1}$ are surjective  for 
$k\ge 2$, which implies surjectivity of $R_k\ot T\ra R_{k+1}$. \qed

\subsection{Their normal spaces}\label{minunormal}

An interesting property of $G$-minuscule varieties is that the
irreducibility of the tangent space propagates to the
irreducibility of all normal spaces. Indeed,
the normal spaces and fundamental forms
of the minuscule varieties are as follows:

\begin{prop}\label{tablesherm}
 The tangent space $T$, and the normal spaces $N_j$, 
with $2\le j\le l$,  of the classical irreducible 
minuscule varieties $X$ 
are given by the following table: \medskip  

\begin{center}\begin{tabular}{||c|c|c|c|c||} \hline 
$X$  &  $G(k,n)$  & $G_{Lag}(n,2n)$ & $\SS_{2n}$ & $\QQ^n$\\ 
$G$  & $SL_n$  & $Sp_{2n}$ & $\Spin_{2n}$ & $SO_{n+2}$ \\  
$H$ &   $SL_k\times SL_{n-k}$  & $SL_n$ & $SL_n$ & $SO_n$\\
$T$ & $W_{\o_{k-1}+\o_{k+1}}=E^*\ot Q$ &  $W_{2\o_1}=S^2U$ & 
$W_{\o_ {n-2}}=\Lambda^2U$ & $W_{\o_2}$\\ 
$N_j$ & $W_{\o_{k-j}+\o_{k+j}}=\La j E^*\ot \La j Q$ & 
$W_{2\o_j}=S_{{\scriptstyle 2\ldots 2}}U$ & $W_{\o_ {n-2j}}=\Lambda^{2j}U$
& $\CC$ \\ 
$l$ & $\min (k,n-k)$ & $n$ & $[\frac{n}{2}]$ &  $2$ \\ \hline   
\end{tabular}\end{center}

\medskip For the two exceptional irreducible minuscule varieties, we have
the following table:\medskip

\begin{center}\begin{tabular}{||c|c|c|c|c|c||} \hline 
 $X$  &  $G$  &  $H$ &  $T$  &  $N_2$ & $N_3$  \\ \hline 
$\BO\pp 2$ &  $E_6$ & $\Spin_{10}$ & $W_{\o_ 2}$ &  $W_{\o_ 6}$ &  $0$  
\\  \hline
$G_{\o}(\BO^3,\BO^6)$  &  $E_7$ & $E_6$ & $W_{\o_6}$ & $W_{\o_1}$ & $\CC$
\\ \hline \end{tabular}\end{center}\end{prop}

The fundamental forms may be described explicitly as follows:
\smallskip

For a non-degenerate quadric $\QQ^n$,   the 
second fundamental form is 
a nondegenerate quadratic form with base locus a smooth quadric 
$\QQ^{n-2}$. 

For the respective cases $G(k,v), G_{\o}(k,V), \BS ,
G_{\o}(\BO^3,\BO^6)$,   $T$ is a (subset) of a matrix space, 
respectively $T= E^*\ot Q, S^2E^*, \La 2 E^*,\JO$. In all but
$\BS$, the last fundamental form is the set of maximal minors
(the determinant  for $ S^2E^*$ and $\JO$), and the lower 
fundamental forms are just the successive Jacobian ideals.
For $\BS$, the last form is the Pfaffian (since the determinant
is a square) and the other forms are the successive Jacobian
ideals, which are the Pfaffians of the minors centered about
the diagonal.

For the $\BO\pp 2$ case, let $V=\bcc {10}$. Then 
$T=S_+(V)$ is a half-spin representation, and
$N_2$ is  the vector representation $V$. The half-spin representations
$S_+$ and $S_-$ can be constructed as the even and odd parts of the 
exterior algebra of a null $5$-plane $E$ in $V$
(see  e.g. \cite{harvey}): 
$S_+, S_-$ are dual to one another, the wedge product
giving a perfect pairing $S_+\otimes S_-
\lra \Lambda^5E = \CC$. Moreover, the full 
exterior algebra of $E$ is a module over the Clifford algebra of $V$. 
If $F$ is a complementary null $5$-plane of $E$, then $E$ acts on $S_+$
by  exterior multiplication, $F$ by interior multiplication, and this
action of $V = E\op F$ extends to the whole Clifford algebra. 
In particular, there is
a natural map from $V$ to $End(S_-,S_+)\simeq S_+\ot S_+$.  
The transpose of the symmetric part of this morphism is the 
second fundamental form.

 Alternatively, identifying $S_+(V)=\BO\oplus\BO$ (see \cite{harvey})
with octonionic coordinates $u,v$, we have $\ii =\{
u\overline u, u\overline v, v\overline v\}$ where,
considering $\BO$ as an eight dimensional vector space
over $\CC$, the middle
equation is eight quadrics.

 \begin{rema}
 Note that in all cases, the {\em only} $H$-orbit closures 
in $\ppp T_1$ are the secant varieties. This actually characterizes
the minuscule varieties, see \cite{LM1}. A special case of this
phenomenon is observed in \cite{rui}.
Note that this property also allows one to easily classify the
$G$-orbits in $\t (X)$ when $X$ is minuscule. See
\cite{LM2} for examples.\end{rema}

\begin{coro}
Let $X$ be a  minuscule variety, and let  $x\in X$.
Then
$$\tbase |\FF^k_{X,x}|  =  \sigma_{k-1}(\tbase |\FF^2_{X,x}|).$$
Moreover, $|\FF^k_{X,x}| = I_k (\tbase |\FF^k_{X,x}|)$.
\end{coro}

\noindent {\em Proof of the corollary}. 
Immediate from our
explicit descriptions of the fundamental forms. \qed

\medskip\noindent {\em Proof of the proposition}. 
For each of these varieties, and each integer $j$, we   check that 
there is a unique irreducible $H$-module which is a component of both
$S^jT$ and  of the restriction $Res^G_HV_{\o_ i}$.  
Then  $N_j$ must be this $H$-module.  

For an ordinary Grassmannian $G(k,n)=G(k,V)$,   $T = E^*\ot Q$,
where $E$ is the tautological subbundle and $Q=V/E$ the
quotient bundle.  Its symmetric
powers are given by the Cauchy formula (\cite{mac}, p. 33)
$$
S^jT   = \bigoplus_{|\l| = j}S_{\l}E^*\ot S_{\l}Q,
$$ 
the
sum is over all partitions $\l$
with  the sum of its  parts $|\l|$   equal to
$j$. We have 
$$
Res^G_H\Lambda^ k V  =  \Lambda^k(E\oplus Q)   = 
\bigoplus_{h\geq 0}\Lambda^hE^*\ot\Lambda^hQ =  
\bigoplus_{h\geq 0}W_{\o_{k-h}+\o_{k+h}} 
$$
since $\trank (E)=k$. 
The only common component of these two decompositions is  
$\Lambda^jE^*\ot\Lambda^ jQ = W_{\o_ {k-j}+\o_ {k+j}}$.  
The case of Lagrangian Grassmannians is  similar. Here   $Q\simeq
E^*$,  $T = S^2E^*$  and we   use the formula (\cite{mac}, p. 45)
$$
S^jT   = \bigoplus_{|\l| = j}S_{2\l}E^*.
$$
We compute the decomposition 
$$
Res^G_HV_{\o_ n}  = 
\bigoplus_{h\geq 0}S_{\underbrace{{\scriptstyle 2\ldots 2}}_{h}}E^*
  =  \bigoplus_{h\geq 0}W_{2\o_h},
$$ 
and the conclusion  follows as above. 
On spinor varieties,  $Q\simeq
E^*$ again, $T = \Lambda^2E^*$ and we  use the formula (\cite{mac}, p. 46)
$$
S^jT  = 
\bigoplus_{|\l| = j}S_{\l(2)}E^*,
$$
where 
if $\l  =  (\l_1,\ldots ,\l_m)$, then   $\l(2) =  
(\l_1,\l_1,\ldots ,\l_m,\l_m)$ . 
Finally, the case of quadrics is immediate since they are
hypersurfaces.  

\medskip For exceptional minuscule varieties the same argument goes 
through, except that we use the LiE package \cite{lie},
or Littelmann paths, instead of the above classical decomposition 
formulas.\qed

\subsection{Algebraic structures induced by 
infinitesimal geometry}\label{algstrs}

We remark on some consequences
of the strict prolongation property
for minuscule varieties.

\begin{prop}\label{algnorm}
 Let $X^n\subset\pp\na$ be a  variety 
such that strict prolongation holds at $x\in X$. Let $N_j=N_{j,x}X$.
  Then there are natural maps
$$N_i^*\ot N_j^*\ra N_{i+j}^*.$$ \end{prop}

\proof The maps are the restrictions of the
symmetrization  maps $S^iT^*\ot S^jT^*\ra
S^{i+j}T^*$ and the image is assured to lie in $N^*_{i+j}$ by
the strict prolongation property.\qed

\begin{coro} Let $X = G/P\subset \PP V$ be a minuscule variety. 
Let $H$ be the semi-simple part of $P$.
Then there is a natural structure of a graded $H$-algebra on $V$. 
\end{coro}

In the case of Grassmannians, the algebra structure on $\La k V$
is given by the multiplication of minors. For Lagrangian
Grassmannians, it is related to the multiplication of Pfaffians. 
The exceptional cases are related to certain exceptional algebraic 
structures introduced by Freudenthal, that we   meet again in
\S 6. For example, consider the minuscule variety of $E_7$: this is a
$27$-dimensional subvariety of the  projectivization of the minimal
representation $V$ of
$E_7$, whose dimension is $56$. As an $H$-module, we have 
  $$V=V_0\op V_1\op V_2\op V_3=\CC\op\JO\op\JO ^*\op\CC,$$
where $\JO$ denotes the exceptional Jordan algebra of $3\times 3$
$\OO$-Hermitian matrices. The group $H=E_6$ is  realized as the
subgroup of $GL(\JO)$ preserving the cubic form defined by the
determinant. Its polarization defines the map $V_1\ot V_1\ra V_2$. 
The map $V_1\ot V_2\ra V_3$ is just the evaluation. 
\medskip

Another consequence of the strict prolongation property
at a point of any variety is
the appearance of Koszul complexes:

\begin{coro}\label{koszul} Let $X^n\subset\pp\na$ be a  variety 
such that strict prolongation holds at $x\in X$. Let
$N_j=N_{j,x}X$.  Then
there is a Koszul complex:
$$
\cdots\lra N_{j-1}^* \ot \Lambda^{k+1}T^*\lra N_j^* \ot \Lambda^kT^*
\lra N_{j+1}^* \ot \Lambda^{k-1}T^*\lra\cdots 
 $$
induced by the maps $T^* \ot N^*_j\ra N^*_{j+1}$ (recall that
$T^*=   N^*_1$).\end{coro}

If $N_j^*$ is replaced by the space of sections $\Gamma
(X,\cO_X(j))$ for a subvariety $X\subset\PP  T $, the homology of the 
corresponding  Koszul complexes compute the syzygies of $X$ \cite{gre}.

\smallskip
For a classical minuscule variety $X$, there is a strange 
relation between the complexes constructed from their normal 
spaces, and the Koszul complexes computing the syzygies of another 
minuscule variety $Z$. Indeed, we obtain this second family of complexes 
from the first,   by a natural
involution on the set of highest weights of irreducible $L$-modules. 

For $L = GL_n$, this involution is defined in the following way:
to the Schur power $S_{\l}$ we associate $S_{\l^*}$, where $\l^*$ is
the conjugate partition of $\l$, obtained by  symmetry along the 
main diagonal of its diagram (which actually defines a bijection 
between partitions inscribed in a $k\times (n-k)$ rectangle, and 
partitions inscribed in a $(n-k)\times k$ rectangle). The complex
associated to our example above is therefore 
$$\cdots \lra 
S^j E\ot S^j F\ot\La k (E\ot F)\lra
S^{j+1}E\ot S^{j+1}F\ot\Lambda^{k-1}(E\ot F)\lra\cdots$$

In small degrees, 
for $X=G(k,n),\;G_{Lag}(n,2n),\;\SS_n$, 
we obtain the Koszul complexes associated to 
$Z = \PP^{n-k-1}\times \PP^{k-1},\;G(2,n),\;v_2(\PP^{n-1})$ respectively. 
(Note that $Y_1 = \PP^{k-1}\times
\PP^{n-k-1},\;v_2(\PP^{n-1}),\;G(2,n)$ respectively.)

\section{Linear spaces on homogeneous varieties}
In this section we explictly describe the lines through a point
of a homogeneous variety $X=G/P\subset\BP V$, the Fano variety
parametrizing all lines on $X$  (in \S\ref{lines}), 
as well as the Fano varieties
parametrizing higher dimensional linear spaces on $X$  in \S\ref{higherdims}.
These parameter spaces all come equipped with natural
embeddings to projective spaces whose associated
vector spaces are $G$-modules as described in \S\ref{embeddings}.
We give an amusing recipe for recovering the Dynkin diagram of
$G$ from   second fundamental form data in \S\ref{ddviaii}.
We begin, in \S\ref{titsfibrations} with a construction due to
Tits that is essential to our work.

\subsection{Tits fibrations}\label{titsfibrations}
Let $G$ be a simple Lie group,
let $S,S'$ be two subsets of the sets of simple roots of
$G$. Consider the diagram 
$$\begin{array}{rcccl} & & G/P_{S\cup S'} & & \\ &
{\scriptstyle \pi} \swarrow & &
\searrow {\scriptstyle \pi '} & \\
X= G/P_S & & & & X'= G/P_{S'}
\end{array}$$
Let $x'\in X'$ and consider $Y:= \pi (\pi'\inv (x'))\subset X$. Then $X$ is
covered by such varieties $Y$. Tits shows in \cite{tits} that $Y= H/Q$ where
${\mathcal D}(H)={\mathcal D}(G)/(S\backslash S')$, and 
$Q\subset H$ is the parabolic subgroup corresponding to $S'\backslash S$.
He  calls such subvarieties $Y$ of $X$, $L$-subvarieties, and $Y$ the 
{\em shadow} of $x'$. 

\begin{exam} For $X= D_n/P_3$ and $X'= D_n/P_{n}$, we read off the diagram
below that $Y= G(3,n)$.\smallskip

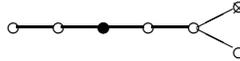
\setlength{\unitlength}{3mm}
\begin{picture}(10,2.5)(-18,-1)
\multiput(0,0)(2,0){5}{$\circ$}
\multiput(0.45,.35)(2,0){4}{\line(1,0){1.65}}
\put(8.45,.25){\line(2,-1){1.6}}
\put(8.45,.35){\line(2,1){1.6}}
\put(4,0){$\bullet$}
\put(10,.9){$\circ$}\put(9.95,.95){${\scriptstyle \times}$}
\put(10,-1.1){$\circ$}\end{picture} \end{exam} \medskip

\subsection{Lines}\label{lines}

Let $\cD=\cD (G)$ be the Dynkin diagram of a complex simple Lie group
$G$. We identify the nodes of
$\cD$ with the set of simple roots with respect to a choice
of maximal torus $T$ and Borel subgroup $B$ which we fix
once and for all.
 Let $\a\in \cD$. Let  $N(\a )=\{\b\in\cD  \mid (\a,\b)< 0\}$ 
denote the {\it  neighbors} of $\a$, the simple roots connected
to $\a$ by an edge in $\cD$. 

\smallskip
Proposition 2.5 implies the following (with the same 
notations):

\begin{coro}\label{wtits}
Let $X=G/P$, with $G$ simple and $P$ a maximal parabolic subgroup. 
Let $Y_1\subset \ppp T_1$ be the closed orbit. 
Then $Y_1$ is isomorphic to the shadow of a point $x\in X$ 
on the space 
$X'=   G/P'$ of $G$-lines in $X$.\end{coro}

Here $P'$ is the parabolic subgroup of $G$ defined by the neighbors
of the root defining $P$ . 

We will see that if we consider the minimal 
homogeneous embedding $X\subset \PP V$, these $G$-Tits lines are
linearly embedded. We first need to recall a few basic facts on 
the Picard group of a rational homogeneous space.

It a classical fact, due to Chevalley, that $Pic(G/B)=H^2(G/B,\ZZ)=P$, the 
weight lattice (\cite{chev}, Expos\'e 15). More generally,  $Pic(G/P_S)=H^2(G/P_S,\ZZ)=P(S)$, 
the sublattice generated by the fundamental weights $\o_i$ dual to the roots
$\a_i\in S$ (or rather to the corresponding coroots).  Dually, $H_2(G/P_S,\ZZ)\simeq \check R(S)$, 
the lattice generated by the coroots $\check\a_i$ to the roots $\a_i\in S$, with the 
obvious pairing with $P(S)$. Each class $\check\a_j$ can be realized geometrically 
by considering the double fibration 
$$\begin{array}{ccccc}
 & & G/P_{S \cup N(\a_j)} & & \\
 & \swarrow & & \searrow & \\
G/P_S & & & & G/P_{S\backslash \a_j \cup N(\a_j)}
\end{array}$$
Indeed, the shadow on $G/P_S$ of a point in $G/P_{S\backslash \a_j \cup N(\a_j)}$ is a rational 
curve, on which a line bundle $L_{\l}$ defined by a weight $\l\in P(S)$ has degree 
$\langle l,\check\a_j\rangle$ (see \cite{dema}, Lemme 2 p. 58). 

In particular, suppose that $G/P_S$ is embedded in some $\PP V_{\l}$ by a very ample
line bundle $L_{\l}$, where $\l=
\sum_{i\in S}l_i\o_i$, and contains a line of $\PP V_{\l}$ whose 
homology class is $\check\b=\sum_{i\in S}m_i\check\a_i$. Then 
$\sum_{i\in  S}l_im_i=1$, which implies that $\check\b=\check\a_j$ for some
$j\in S$ with $l_j=1$.  Moreover, the variety $F^j_1(X)$ of these 
$\a_j$-lines on $X$ is independent of the $\l$ with $l_j=1$ chosen.
Note that is contains $G/P_{S\backslash \a_j \cup N(\a_j)}$.

\begin{theo}\label{line}
Let  $S\subseteq\cD$, consider $X=G/P_S$ in its minimal homogeneous 
embedding. Then 

1. $F_1(X)=\coprod_{j\in S}F_1^j(X)$, where $F_1^j(X)$ is 
the space of lines of class $\check\a_j\in H_2(G/P_S,\ZZ)$.

2. If $\a_j$ is not an exposed short root, then $F_1^j(X)=
G/P_{S\backslash \a_j\cup N(\a_j)}$.

3. If $\a_j$ is an exposed short root, then $F_1^j(X)$ is the union of 
two $G$-orbits, an open orbit and  its boundary $G/P_{S\backslash\a_j\cup N(\a_j)}$.
\end{theo}

Assertions 1. and 2. are rephrasings of results in \cite{coco}, which were 
published just after the first version of this paper was written (but note 
that Cohen and Cooperstein work over an arbitrary field). Assertion 
3. and its proof below were communicated to us by an anonomous referee
(our original proof contained some case by case arguments).

In \S 6 we   give explicit descriptions of the open orbits 
of 3. for each short root.

\proof The argument,
  proceeds in three
 steps: first we give a criterion for identifying distinct
orbits in $F^j_1(X)$; up to the action of $W_S$,
the subgroup of the Weyl group $W$ generated by
the simple reflections $s_i$, $i\notin S$, there    is a unique $T$-fixed
point in each $G$-orbit passing through a base point $x\in X$, where $T$ denotes the maximal torus in $G$.
We then show that there is a unique tangent direction in $T_{j1}$ correspoinding to a line up to
$W_S$-equivalence, and finally, if $\a_j$ is an exposed short root
with $-(\a,\a_j)=2$ (resp. $3$), then there is a unique vector in $T_{j2}$
(resp. $T_{j3}$) modulo $W_S$-equivalence corresponding
to a $T$-fixed line. Finally we show that the orbit of the $T_{j2}$ (resp. $T_{j3}$) line is not
closed.

We first observe that every $G$-orbit in $F^j_1(X)$ contains
a $T$-fixed line. Indeed, if $l$ is an $\a_j$-line, by homogeneity
we can suppose it contains the base point $x=[P_S]\in X$. If $y\in l$
is another point, $y$ lies in a unique Bruhat cell $B(w)=
\{ wxh\inv \mid h\in B\}$, $w\in W$.
Then there exists $h\in B$ such that $h(y)=w(x)$, and since $h(x)=x$,
the line $l'=\overline{xw(x)}$ is  $T$-stable  in   $Gl$. 

By $T$-stability, the tangent direction to $l'$ at $x$ must be a 
$\fg_{-\a}$ for some $\a\in\Delta_X$.
The orbit through $x$  of the $SL_2$ corresponding to $\a$ is $l$,
and thus $\langle \l,\check\a \rangle =1$ because $\lambda$ is also the highest
weight of the $SL_2$-module $\hat l\simeq\CC^2$.

Write $\a=p\a_j +\g$ with $m_j(\g)=0$. We have
$$1=\langle \l, \check\a \rangle  = p\frac{(\a_j,\a_j)}{(\a,\a)}
+ \Sigma_{i\ne j} l_i\langle \o_i, \check\a\rangle .$$
Since the equality holds for any choices of coefficients
 $l_i$, we must have 
$\langle \o_i, \check\a\rangle =0$ for all $i\in S\backslash j$.
Thus  $(\a,\a)=p(\a_j,\a_j)$ and we have two cases. 

If $\a_j$ is not an exposed short root, then $p=1$ and by Lemma 
\ref{slem} below, $\a$ is $W_S$-conjugate to $\a_j$. There is
therefore a unique $G$-orbit in $F^j_1(X)$. 

If $\a_j$ is an exposed short root, then either $p=1$ and $\a$ is
short, or $p>1$ and $\a$ is long. Then Lemma \ref{slem} and Lemma
\ref{lemp>1} below imply that there are at most two $G$-orbits in $F_1^j(X)$.

\begin{lemm}\label{slem}
Let $\b\in \Delta_X$ have the same length as $\a_i$, and $m_i(\b)=1$
Then there exists $w\in W_S$ with $w\a_i=\b$. 
\end{lemm}

\proof We use induction on the height
of $\b$, i.e., the sum of its coefficients in its decomposition 
on simple roots. Suppose that $\b\ne\a_i$, 
and  $(\b,\a_j)\le 0$ for all $j\ne i$. Then
$(\b,\a_i)\ge (\b,\b)>0$, thus $n(\b,\a_i)=1$ and  $s_i(\b)=\g$, where
 $s_i(\b)= \b - \langle \b  , \check \a_i\rangle \a _i$.
Therefore, since $\b,\g$ and $\a_i$ have the same length, we get 
$-1 = n(\g,\a_i)=n(\a_i,\g)$, hence $n(\b,\g)=1$ and $(\b,\g)>0$.

We can therefore let $k\ne i$ such that $(\b,\a_k)>0$. Then the 
root $s_k(\b)$ verifies the same assumptions as $\b$, but its height
is smaller, and we conclude the proof by induction. \qed

\begin{lemm}\label{flem} $\a_j$ is an exposed short root
iff there exists $\a\in\Delta_+\backslash S$
such that $|n(\a, \a_j)|>1$.\end{lemm}

\proof If $\a_j$ is not an exposed short root, $n(\a,\a_j)
\neq 0$ implies that $\a$ and $\a_j$ have the 
same length, hence $|n(\a, \a_j)|=
|n(\a_j, \a)|\le 1$. 
Now say $\a_j$ is  an exposed short root. Then there exists a long 
root $\a$ supported outside $S$ such that $(\a,\a_j)\neq 0$. Thus 
$|n(\a,\a_j)|>|n(\a_j,\a)|=1$.\qed

\begin{lemm}\label{lemp>1}
If $\a_j$ is an exposed short root, then any pair of long roots of the form $p\a_j+\g$, where $\g$ is
supported outside $S$, are conjugate in $W_S$.
\end{lemm}

\proof If $p=3$, we are in the $G_2$ case which is clear. 
So assume $p=2$. Let $\d= 2\a_j +\g$, $\g$ supported outside $S$.

Assume $\d$ is such that $(\d, \a_i)\leq 0$ for each $i\neq j$ with 
$m_i(\d)\neq 0$. Then $(\d, \a_j)\ge (\d,\d)>0$.
Since $\a_j$ is short, we have $n(\d,\a_j) =2$,
so $s_j(\d)=\g$ and $\g$ is a long root. Also, $n(\a_j,\d)=1$, 
so that $n(\d,\a_i)=0$ for all $i\neq
j$  with $m_j(\d )\neq 0$. Thus $\d$ is the highest root with support
in the subdiagram $supp(\d )$. But $n(\d,\a_j)=2$,
therefore this subdiagram must be of type $C_{r+2}$ for some $r\ge 0$. 
Since $j$ is an end of it, it is uniquely determined, and $\d$ is a
uniquely determined root $\d_0$. 

Thus if $\d\ne\d_0$, there 
exists $i\neq j$, $\a_i\in \cD\backslash S$, with $(\d,\a_i)>0$. And
we can then proceed by induction on the root $s_i(\d)$, whose height 
is smaller than that of $\d$. \qed

\medskip
Finally we show that when $\a_j$ is an exposed short root, 
the orbit corresponding to a long root $\a$ is not closed. 
The case of $G_2$ is easy (see \S 6.1), so assume $p=2$
and, using the notation of lemma \ref{flem}, take $\a=\d_0$.
The corresponding line is $l_{\a}=\overline{v_{\l}v_{\l-\a}}$, where 
$v_{\l}$ is a highest weight vector in $V_{\l}$ and $v_{\l-\a}=X_{-\a}v_{\l}$, 
$X_{-\a}\in \fg_{-\a}$. Let us compute the tangent space $T_{\a}$ 
at $l_{\a}$ of its $G$-orbit. For $X\in\fg$, we have
$$X(v_{\l}\we v_{\l-\a})=Xv_{\l}\we v_{\l-\a}+v_{\l}\we
XX_{-\a}v_{\l}.$$
If $X\in\fg_{-\b}$, $\a\ne\b\in\Delta_X$, then $Xv_{\l}\we v_{\l-\a}$,
and therefore also $X(v_{\l}\we v_{\l-\a})$, are non zero. If
$X\in\fp_S$, then $Xv_{\l}=0$, hence $X(v_{\l}\we v_{\l-\a})=
v_{\l}\we [X,X_{-\a}]v_{\l}$. If $X$ is a root vector, this is non 
zero if and only if $X\in\fg_{\a-\b}$, with $\b\in\Delta_X$. 
This implies that, as a $T$-module,  
$$T_{\a}=\bigoplus_{\b\in\Gamma}\fg_{-\b}, \quad \Gamma = (\Delta_X-\{\a\})\cup
\{\b-\a\in\Delta, \b\in\Delta_X\}.$$
(Note that since $\a$ is a long root, 
we also have $\Gamma = (\Delta_X-\{\a\})\cup s_{\a}(\Delta_X-\{\a\})$.) 

To prove that the $G$-orbit of $l_{\a}$ is not closed, we need to
prove that its stabilizer cannot contain a Borel subgroup, and for
this it is enough to exhibit a root $\b\in\Gamma$ such that $-\b$ also 
belongs to $\Gamma$. But this is easy: indeed, recall that the
subdiagram of $\cD$ supporting $\a$ is of type $C_{r+2}$, and that 
$\a$ is the corresponding highest root, that is $\a=2\a_j+\cdots +
2\a_{r+j}+\a_{r+j+1}$ for a suitable numbering of the simple roots. 
Then we can take  $\b=\a_j+\cdots+\a_{r+j}+\a_{r+j+1}\in\Delta_X$, 
since $-\b+\a=\a_j+\cdots+\a_{r+j}\in\Delta_+$. \qed

\medskip 
\begin{rema} This proof gives a formula for the
dimension of the open orbit of lines in the case of an exposed short
root, namely the cardinality of $\Gamma = (\Delta_X-\{\a\})\cup 
s_{\a}(\Delta_X-\{\a\})$. \end{rema}

\medskip
On the infinitesimal level, $C_x\subset\ppp T_xX$, the set
of tangent directions to lines on $X$ passing through $x$, is a union of disjoint
varieties, one component $C_x^\a$ for each possible class $\a$ of lines.
The proof of the preceeding theorem implies the following:

\begin{theo}\label{conex}  Let $G$ be a complex simple Lie group,
let $S$ be a subset of the simple roots.  Let $\a\in S$.
 Let $\cD (H)$ be the components of $(\overline{\cD (G)\backslash S})
\backslash\a$ containing
an element of $N(\a )$, where by $\overline{\cD (G)\backslash S}$ we
mean $\cD (G)\backslash S$ plus any nodes of $S$ attached
to a node of $\cD (G)\backslash S$. Let $C_x^\a\subset\BP T_xX$
denote the class of $\a$-lines through $x$.

1. If $\a$ is not an exposed short root, then
$C_x^{\a}\simeq H/P_{N(\a )}$.

2. If $\a$ is  an exposed short root, then $C_x^{\a}$ is a union of an 
open $P_S$-orbit and its boundary $H/P_{N(\a )}$.\end{theo}

\medskip
With the notations of \S 2, the closed $P_S$-orbit in $C_x^{\a_j}$
is $Y_1^j$, and the open orbit is $P_SY_2^j$.
The cases of  $C_x^\a$ for exposed short roots are described
explicitly on a case by case basis in \S 6.

\subsection{Linear spaces of higher dimension}\label{higherdims}

A $k$-plane in $X$ must come from a linear $\pp{k-1}$ in some $C_x^\a$.
We call such a $\pp k$, of class $\a$, and let $F^{\a}_k(X)$ denote
the variety parametrizing the $\a$-class $\pp k$'s on $X$.
$F^{\a}_k(X)$ may have several components.
The space of $\pp{k}$'s
in $X$, $F_k(X)$,  is the disjoint union of the $F_k^{\a}(X)$'s.

If $\a=\a_j\in S$ is not an exposed short root, it follows from Theorem
4.3 that the projection $G/P_S\rightarrow G/P_{S\backslash j}$ is constant 
on each $\a$-line, hence it is also constant on each $\a$-class $\pp k$. 
It follows that the space of $\a$-class $\pp k$'s is a fibration over 
$G/P_{S\backslash j}$, and to determine the fiber, we can restrict to the 
subdiagram of $\cD$ consisting in the connected component of $\a_j$ 
in $\cD\backslash (S\backslash j)$. In particular, 
we are reduced to the case of a maximal
parabolic subgroup, corresponding to a non-short root.  

Then we know that $C_x=Y_1$ is a minuscule variety, so it is again 
a homogeneous space of type $H/Q$ with $Q$ a maximal parabolic
subgroup corresponding to a long root, or possibly a product of such
spaces. We can therefore apply Theorem 4.3 to $Y_1$ to describe
its lines, which   gives $\PP^2$'s in the original space, and so
on. The conclusion is that, not only $\PP^1$'s, but all linear spaces can 
be described in terms of Tits' geometries.

\begin{theo}\label{kplanes1}
Let $G$ be a simple group and
let $X=G/P_S\subset\ppp V$ be  a rational homogeneous
variety in its minimal homogeneous embedding.

If $\a\in S$ is not an exposed short root, then for all $k$,
$F^{\a}_k(X)$ is the disjoint union of homogeneous varieties 
$G/P_{\Sigma \b_j}$
where $\{\b_j\}\subset\Delta_+$ is a minimal set of positive roots such that
the component of ${\mathcal D}(G)\backslash \{\b_j\}$ containing $\a$ is
isomorphic to ${\mathcal D}(A_k)$, intersects $S$ only in $\a$, and 
$\a$ is an extremal node of this component.
\end{theo}

\begin{coro} Let $G$ be a simple Lie group, let $S\subset \cD (G)$,
let $\a\in S$ with $\a$ not exposed short.   
Let $X=G/P_S$ be the corresponding homogeneous
variety in a homogeneous embedding such that there are
$\check\a$-lines. Suppose that  the longest
of the chains of type $A$
in $\cD (G)$ beginning at $\a$  and containing no other element of $S$, 
is isomorphic to  ${\mathcal D}(A_n)$.
Then  the largest  linear space of class $\a$ on $X$ is a $\pp n$.  \end{coro}

\begin{exam} Consider the case of $G/B\subset\ppp V_{\l}$,
 $\l= \o_1+\cdots + \o_r$, the sum of the fundamental weights.
There are no unexposed short roots, so $C_x$ is the union of $r$
points, $F_1(X)=\coprod_{j} G/P_{\cD\backslash\a_j}$
and $F_k(X)=\emptyset$ for $k\geq 2$.\end{exam}

\begin{exam} In the case of $D_n/P_n$, we have a unique family of 
lines, parametrized by isotropic subspaces of dimension $n-2$, $D_n/P_{n-2}$
and a unique family of two planes parametrized by the $Q$-isotropic flag
variety $D_n/P_{n-3,n-1}$. There are two families of $\pp 3$'s,
namely $D_n/P_{n-3}$ and $D_n/P_{n-4,n-1}$. For $4\leq k\leq n-1$
there is a unique family of $\pp k$'s, namely $D_n/P_{n-k-1, n-1}$. 
\end{exam}

\begin{exam} The largest linear space on $E_n/P_1$ is a $\pp{n-1}$,
via the chain terminating with $\al n$, so $E_n/P_1$ is  
maximally uniruled
by $\pp{n-1}$'s and there is a second chain terminating with
$\al 2$, so $E_n/P_1$ is also maximally uniruled by $\pp 4$'s. (The 
unirulings by the $\pp 4$'s are maximal in the sense that none
of the $\pp 4$'s of the uniruling are contained in any  
$\pp{5}$ on $E_n/P_1$.)  The varieties parametrizing these rulings are
respectively $E_n/P_2$ and $E_n/P_n$. \end{exam}

\medskip

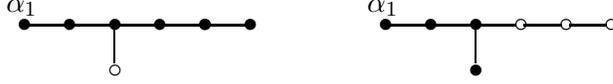
\setlength{\unitlength}{3mm}
\begin{picture}(20,3.5)(-10,-1.8)
\multiput(0,0)(2,0){6}{$\bullet$}
\multiput(0.45,.3)(2,0){5}{\line(1,0){1.65}}
\put(8,0){$\bullet$}
\put(4,-2){$\circ$}
\put(4.3,-1.4){\line(0,1){1.5}}
\put(-.5,.9){$\a_1$}

\multiput(16,0)(2,0){3}{$\bullet$}
\multiput(22,0)(2,0){3}{$\circ$}
\multiput(16.45,.3)(2,0){5}{\line(1,0){1.65}}
\put(20,-2){$\bullet$}
\put(20.3,-1.4){\line(0,1){1.5}}
\put(15.5,.9){$\a_1$}
\end{picture}
 
\bigskip
Now we address the case of exposed short roots. First note
that if $X=G/P_{\a}= B_n/P_n, G_2/P_1$ or $C_n/P_1$,
then the space of $\pp k$'s on $X$ is  $\tilde G$-homogeneous, where ${\mathcal D}(G)$ is the
fold of ${\mathcal D}(\tilde G)$, as in these cases $G/P\simeq
\tilde G/\tilde P$. In general, we have: 

\begin{theo}
If $\a\in S$ is an exposed short root, then for all $k$,
$F^{\a}_k(X)$ consists
 of a finite number of $G$-orbits (at least two).
\end{theo}

If $\a\in S$ is an exposed short root, $F^{\a}_k(X)$ can be deduced
from   $F_k(G/P_{\a})$. In each of these cases we determine
the unextendable linear spaces through a point explicitly. By further
calculation, one can deduce all linear spaces through a point and
prove the theorem.

\subsection{Natural embeddings of linear spaces}\label{embeddings}

  For $X=G/P_S\subset\ppp V$,  we have
$F_{k-1}(X)\subset G(k,V)\subset \ppp\La kV$.
Thus the connected components of the $F_{k-1}(G/P_S)$'s
come naturally embedded in  some
 irreducible component of $\La k V$  with highest weight
$\lambda$ supported on the weights dual
to the roots appearing in the (unique) closed orbit $G/P_{S'}$
consisting of $G$-homogeneous $\pp k$'s in the component (i.e. $\pp k$'s that are
$L$-varieties in the sense of Tits). While
  $S'$ can be determined pictorially, the multiplicities
of the weights  in general cannot. We now determine the multiplicities in several cases,
in particular the elementary representations defined below.

Fix an end of the Dynkin diagram of $G$,
and label the end node $\a_1$. Following \cite{dynkin},
define the {\em branch} of $\a_1$,
 ${\mathcal B}(\a_1)$,   as the largest 
  chain in ${\mathcal D}(G)$ containing $\a_1$ that is isomorphic to 
${\mathcal D}(A_p)$ or ${\mathcal D}(C_p)$,  such that
  no node    in ${\mathcal B}(\a_1)$ before the last   has valence
three.  We say such a branch has {\em length} $p$.
We label the roots on  ${\mathcal B}(\a_1)$ as $\al 1\hd\a_p$ and
  denote  the fundamental representation corresponding to 
 $\omega_1$ by $V=V_{end}=V_{\lo 1}$. Such an irreducible
 representation is called an {\em elementary representation} in
 \cite{dynkin}. 
  
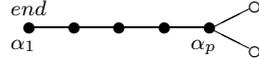
\begin{center}
\setlength{\unitlength}{3mm} 
\begin{picture}(10,3)(0.5,-1)
\multiput(0,0)(2,0){5}{$\bullet$}
\multiput(0.45,.35)(2,0){4}{\line(1,0){1.65}}
\put(8.48,.25){\line(2,-1){1.6}}
\put(8.48,.35){\line(2,1){1.6}}
\put(-.5,.9){$\scriptstyle{end}$}
\put(-.5,-.6){$\scriptstyle{\a_1}$}
\put(7.5,-.6){$\scriptstyle{\a_p}$}
\put(10,.95){$\circ$}
\put(10,-1.05){$\circ$}\end{picture}
\end{center}

\medskip

\medskip
 The following result is evidently due
to Cartan, a proof can be found in \cite{dynkin} except for
the   \lq moreover\rq\ assertions   which may be verified on a case
by case basis.

\begin{prop}  \label{wedge}
With the notations above, 
$V_{\lo k}$ is an irreducible component of $\Lambda^kV$, with
multiplicity one for $2\leq k\leq p$. More precisely, $\o_k$ is the 
unique extremal weight of $\La kV$.

Moreover, $\La {p+1}V$ also has a unique extremal weight which is
 
1.  $2\lo {p+1}$ for a double edge with arrow   pointing away from  $\lo
1$.

2.  $3\lo {p+1}$ for a triple edge with arrow   pointing away from  $\lo
1$.

3. $\lo{p+1}+\lo{p+2}$ if $\lo p$ corresponds to a node of
valence three.
\end{prop}

\medskip\noindent {\em Idea of proof}. 
One simply checks that  among the weights of $V$, there is a maximal 
chain $\mu_1,\ldots ,\mu_{p+1}$ with  $\mu_i = \o_1-(\a_1+\cdots
+\a_{i-1})$. In particular, $\mu_1+\cdots +\mu_k$ is the unique 
maximal weight of $\La k V$ for $1\le k\le p+1$, and it is 
straightforward to check that this weight is as announced in the 
proposition. \qed 

\medskip
It is an easy exercise to prove that
the wedge product of the weight vectors corresponding
to the weights $\m_1,\ldots ,\m_k$ generate a $\PP^{k-1}$ that is 
contained in the closed orbit $X_{end}$. Thus the $G$-submodule
of $\La kV$ which hosts $F_{k-1}(X)$ is precisely the fundamental 
representation $V_{\o_k}$.

It follows for example that in the case of simply-laced groups, 
we can obtain all fundamental representations from the elementary
ones, in a simple geometric way.

\begin{exam} For $E_6$ we have three elementary representations, the 
minimal representation $V_{\o_1}$, its dual
$V_{\o_6}$ and the adjoint 
representation $V_{\o_2}$. 

Start  with $V_{\o_1}$, so $X=E_6/P_1$ is the
Cayley plane.  Then  $\ppp V_{\o_3}$   is 
the ambient space for $F_1(X)$, 
$\ppp V_{\o_4}$   is 
the ambient space for $F_2(X)$
and
$\ppp V_{\o_2}$   is 
the ambient space for $F_5(X)$.

Start with $V_{\o_2}=\fe_6$, so $X=E_6/P_2$ is the adjoint variety.
Then $\ppp V_{\o_4}$   is 
the ambient space for $F_1(X)$    and 
$\ppp V_{\o_3}$   is 
the ambient space for $F_4(X)$. (Note that $E_6/P_1$ is a space
of spinor varieties $D_5/P_5$ on $X$.)
\end{exam}
   
\medskip
A case by case analysis  with LiE \cite{lie} leads to the following 
more precise result:
 
\begin{prop} Notations  as above.

1. If $X_{end}$ is a minuscule variety  and $V$ is not symplectic
or $(D_n,\o_n)\simeq (D_n,\o_{n-1})$,
then for 
$2\leq k\leq p$,  $V_{\lo k}=\La k V$. This is the case for $(G,\o) =
(A_n,\o_1),\; (B_n,\o_1),\;(D_n,\o_1)$ and $(E_6,\o_1)\simeq (E_6,\o_6)$.

2. If $X_{end}$ is  a minuscule variety  and $V$ has a symplectic
form $\Omega$, set $\Lambda^{\langle k\rangle }V = 
\Lambda^{k}V/(\oo\ww\Lambda^{k-2}V)$.
Then $V_{end-k}=\Lambda^{\langle k\rangle }V $. This is the case for
$(G,\o) = (C_n,\o_1)$ and $(E_7,\o_7)$.
 
3.  If $X_{end}$ is an  adjoint variety, so $V=\fg$,   let
$\Lambda^{[2]}\fg = 
\tker [\;,\;]$, where we consider the Lie bracket as
  a  map $[\;,\;] : \Lambda^2\fg\ra\fg$.  
The adjoint varietes corresponding to elementary
representations are those of the exceptional groups: $(G,\o) = 
(G_2,\o_2),\; (F_4,\o_1),\;(E_6,\o_2),\;(E_7,\o_1)$ and $(E_8,\o_8)$.
In each of these cases, except for $G_2$, $V_{end-1}=\Lambda^{[2]}\fg$.
(And afterwards there is always a double bond or node with
triple valence.)

4. If $X_{end}= G_2/P_1$ then $V_{\lo 2}= \La 2 V_{\lo 1}/(V^*_{\lo 1}
\lrcorner\phi)$ where $\phi\in\La 3V_{\lo 1}$ is the defining
three form.

5. If $X_{end}= F_4/P_4$, then $V_{\o_2}=\La 2 V_{\lo 4}/{\mathfrak f_4}$.
\end{prop}

\begin{rema}\label{slad} 
 Consider the case of $A_n$.   The adjoint variety
$X = G/P_{1,n}$  is the flag variety of lines in
hyperplanes in $\pp n$. The space of lines in $X$ is
disconnected: it is the disjoint union of $\FF_{2,n}$ and 
$\FF_{1,n-1}$, the corresponding embeddings of which are {\it not} their 
minimal ones, since 
$$\Lambda^{[2]}\fsl_{n+1} = V_{2\o_1+\o_{n-1}}\oplus V_{\o_2+2\o_n}.$$
\end{rema}

\begin{rema}\label{f4embed}
For $X=F_4/P_4$, $F_5(X)=F_4/P_1$,
but here the variety occurs in its {\it third} Veronese
re-embedding, i.e. $V_{3\o_1}\subset \La 3 V_{\o_4}$.
For $X=G_2/P_1$, $F_2(X)$ occurs in  its second Veronese embedding,
$\ppp V_{2\o_1}$.
\end{rema}

\subsection{Dynkin diagrams via second fundamental forms}\label{ddviaii}

We describe how to recover ${\mathcal D}^*(G)$,  the 
marked Dynkin diagram of
$G$, from the second fundamental form at a point of any $X=G/P$ where
$P$ is maximal and not short:

Fix $x\in X$ and start with a 
marked node, which   corresponds  to $P$.
Say $Y_1$ is the Segre product of  Veronese re-embeddings
of $k$   minuscule varieties.
 Then attach $k$ edges to the node, with
nodes at the end of each edge.
For each factor in the Segre that is minimally embedded,
the edge is simple. If there is a factor  that is a
quadratic (resp. cubic) Veronese, then make
the corresponding edge a double (resp. triple) bond. 
Now compute $\tbase |\fff 2{X,x} |$ of each factor and repeat the
process starting with the node corresponding to the
factor. Continue until arriving at the empty set.
The resulting diagram is ${\mathcal D}(G)$.

A shortcut: if at any point one obtains
an $H/Q$ as a factor where   the marked
Dynkin diagram associated to $H/Q$ is known, one can simply attach the
diagram. In particular, if one arrives at a $\PP^l$, just attach
a copy of ${\mathcal D}(A_l)$.

\begin{exam}  Beginning with $X= E_n/P_1$,
one has 
$$\begin{array}{ll}
\tbase |\fff 2{X}|=\BS_{n-1}, & \tbase |\fff 2{\BS_{n-1}} |= 
G(2,{n-1}), \\ \tbase |\fff 2{G(2,{n-1})}|=
Seg(\pp 1\times\pp {n-2}), & 
\tbase |\fff 2{Seg(\pp 1\times\pp {n-2})}|=
\pp 0\sqcup \pp{n-3}.
\end{array}$$ So the construction is:\medskip

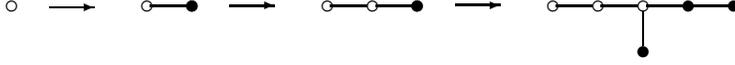
\setlength{\unitlength}{3mm}
\begin{picture}(20,3.5)(-6,-1.8)
\put(0,0){$\circ$}
\put(2,.3){\vector(1,0){2}}
\multiput(6,0)(2,0){2}{$\circ$}
\put(8,0){$\bullet$}
\put(6.45,.35){\line(1,0){1.6}}
\put(10,.35){\vector(1,0){2}}
\multiput(14,0)(2,0){3}{$\circ$}
\put(18,0){$\bullet$}
\multiput(14.45,.35)(2,0){2}{\line(1,0){1.6}}
\put(20,.4){\vector(1,0){2}}
\multiput(24,0)(2,0){3}{$\circ$}
\multiput(30,0)(2,0){2}{$\bullet$}
\put(28,-2){$\bullet$}
\multiput(24.45,.35)(2,0){4}{\line(1,0){1.6}}
\put(28.3,-1.4){\line(0,1){1.5}}
\end{picture}\end{exam}

\section{Classical homogeneous varieties}

In this section  we  present the higher normal spaces of the 
classical homogeneous varieties which are not minuscule, as well
as the base loci of the higher fundamental forms. We give most results
without their proofs, which are computational. 

\subsection{Orthogonal Grassmannians}

Let   $G_o(k,n)$ denote the
orthogonal Grassmannian of null $k$-planes in $V = \CC^n$ 
where $V$ is equipped with a nondegenerate
quadratic form $Q$ and  $2k< n$. It is a
subvariety of the ordinary Grassmannian, and its minimal embedding is 
the Pl\"ucker embedding in $\PP  V_{\o_k}  = \PP (\Lambda^kV)$. 

Introduce the notation
$$
(\Lambda^pE^*\ot\Lambda^pE^*)_+= \left\{ \begin{matrix} 
S^2(\La p E^*) & {\rm if }\ p\ {\rm  is}\ {\rm  even}\\
\La 2(\La p E^*) & {\rm if }\ p\ {\rm  is}\ {\rm odd}\end{matrix}\right .
$$ 
  and
$$
(\Lambda^pE^*\ot\Lambda^pE^*)_-= \left\{ \begin{matrix} 
\La 2(\La p E^*) & {\rm if }\ p\ {\rm  is}\ {\rm even}\\
S^2(\La p E^*) & {\rm if }\ p\ {\rm  is}\ {\rm odd.}\end{matrix} \right .
$$ 

\begin{prop} Let $E^k $ be the   tautological 
vector subbundle  on $G_o(k,n)$,  let $E\upperp\supset E$ denote
its $Q$-orthogonal complement, and let $U^{n-2k}
= E^{\perp}/E$. Then the tangent  space and  normal spaces of $G_o(k,n)$,
as $H= SL(E)\times SO(U)$ modules, are
$$\begin{array}{rcl}
T_1 &  =  & E^*\ot U, \qquad
T_2 \;\; =  \;\; \Lambda^2 E^*, 
\end{array}$$
$$\begin{array}{rcl}
N_2 &  =  & (\Lambda^2 E^*\ot \Lambda^2 U\oplus S^2E^*)\oplus (\Lambda^2
E^*\ot
E^*\ot U)
\oplus (\Lambda^4 E^*\oplus S_{22}E^*).
\end{array}$$
$$\begin{array}{rcl}
N_p  &  =  & \bigoplus_{a>0}(\Lambda^{p-a}E^*\ot\Lambda^pE^*\ot\Lambda^aU) 
\oplus  \\ & & \hspace*{3cm} \oplus (\Lambda^pE^*\ot\Lambda^pE^*)_+\oplus
(\Lambda^{p-1}E^*\ot\Lambda^{p-1}E^*)_-.\end{array}
$$

In particular, the length of the normal graduation is 
$k$ when $k$ is even and the last non zero term is
$N_k\simeq
\La k E^*\ot (\La 0 E^*\ot \La k U +
\hdots + \La k E^*\ot \La 0 U) \oplus S^2(\La k E^*)$.
When $k$ is odd, the length is  $k+1$  and the last non zero term is
$N_{k+1}\simeq\CC$.
\end{prop}
 
\begin{rema}\label{noprolong} 
Note that in contrast to the
case of minuscule varieties, here $N_{3}^*\ne N_2^{*(1)}$.  \end{rema}

\begin{coro} The base locus $\tbase |\fff p{ G_{o}(k,n), E}|$ of the 
$p$-th fundamental form is, for $p$ even, 
$$\ppp \overline{\{e_1\ot u_1+\cdots +e_p\ot u_p\op e_1\we e_2+\cdots 
+e_{p-1}\we e_p, \mid e_j\in E^*,\;u_j\in U\} }, $$
and for $p$ odd,
$$\ppp \overline{\{e_1\ot u_1+\cdots +e_p\ot u_p\op e_1\we e_2+\cdots 
+e_{p-2}\we e_{p-1}, \mid e_j\in E^*,\; \;u_j\in U\} }.$$
\end{coro}

\subsection{Symplectic Grassmannians}

We let  $G_\o(k,2n)=C_n/P_k$ denote
  the Grassmanian of $k$-planes
isotropic for a symplectic form. Its minimal embedding is to 
$V_{\o_k} = \Lambda^{\langle k\rangle }V =
\Lambda^kV/(\Omega\we\Lambda^{k-2}V)$, 
 the $k$-th reduced exterior power of $V=\bcc{2n}$,
where $\Omega$ denotes the symplectic form. 

  A straightforward computation 
shows that   $V_{\o_k}$   
has the following decomposition  as an  $H=SL_k\times
Sp_{2n-2k}$-module:
$$
\Lambda^{\langle k\rangle }V =
\bigoplus_{a,b}\Lambda^bE^*\ot\Lambda^{a+b}E^*\ot\Lambda^{\langle
a\rangle}U.
$$
Note that $U = E^{\perp}/E$ is endowed with a symplectic form induced 
by the symplectic form on $V = \CC^{2n}$.

\begin{prop} Let $E $ be the   tautological 
vector sub bundle  on $G_{\o}(k,2n)$, 
let $E\upperp\supset E$ denote the $\Omega$-orthogonal complement
to $E$ and let $U =
E^{\perp}/E$. Then the tangent  space and normal spaces of $G_{\o}(k,2n)$
are, as
$H$-modules,
$$\begin{array}{rcl}
T_1 &  =  & E^*\ot U, \qquad T_2 \;\;  =  \;\; S^2 E^*, 
\end{array}$$
$$\begin{array}{rcl}
N_2 &  =  & \Lambda^2 E^*\ot \Lambda^{\langle 2\rangle} U\oplus 
S_{21} E^*\ot U\oplus S_{22}E^*, \\
 & & \\
N_p &  =  & \bigoplus_{a+b+c = p}\Lambda^{\langle a\rangle}U\ot S_{
\underbrace{\scriptstyle{2\ldots 2}}_{b-c} 
\underbrace{\scriptstyle{1\ldots 1}}_{a+2c}}E^* \\
&  =  & \bigoplus_{d+e = p}\Lambda^dU\ot S_{
\underbrace{\scriptstyle{2\ldots 2}}_e \underbrace{
\scriptstyle{1\ldots 1}}_d}E^*.
\end{array}$$
In particular, the length of the normal graduation is equal to 
$k$, the last non zero term being $N_k\simeq\Lambda^k(\CC\op U)$.
\end{prop}

\begin{coro} 
$$
\tbase |\fff 2{ G_{\o}(k,2n), E}|= 
\ppp \overline{\{e\ot u\oplus e^2 \mid e\in E^*\backslash \{0\},\;u\in U
\backslash \{0\}\} }. 
$$ 
This base locus contains an open and dense $P$-orbit, 
the boundary of which is the union of the two (disjoint) closed 
$H$-orbits 
$$Y_1\simeq \PP^{k-1}\times\PP^{2n-2k-1}\subset \PP(T_1)\quad
and \quad Y_2\simeq v_2(\PP^{k-1})\subset \PP(T_2).$$
\end{coro}

\proof This  can be seen directly. A line in
$G_{\o}(k,2n)$ through a   point $E$ is given by a $(k-1)$-plane
$H\subset E$, and a 
$(k+1)$-plane $K\supset E$.   $K$ does not need to be isotropic,
each point
of the corresponding line   is generated by $H$ and a vector of $K$, 
and
   is isotropic if and only if  this vector is $\o$-orthogonal to
$H$. The condition on $K$ is
thus $K\subset H^{\perp}$. $K$ is
therefore determined by a line in 
$H^{\perp}/E\simeq U\oplus H^{\perp}/E^{\perp}.$

If $e\in E^*$ is an equation of $H$, the line 
$H^{\perp}/E^{\perp}\simeq (E/H)^*\subset E^*$ is generated by $E$, 
so that a vector in $H^{\perp}/E$ can be written as $u\oplus\l e$, 
where $u\in E$. Our claim follows, the closed orbits $Y_1$ and $Y_2$ 
corresponding to the cases where $u$ or $\l$ is equal to zero. \qed

\begin{coro}
More generally, the base locus of the $p$-th fundamental form is 
$$\tbase |\fff p{ G_{\o}(k,2n), E}|  = 
\ppp \overline{\{e_1\ot u_1+\cdots +e_p\ot u_p\op e_1^2+\cdots 
+e_p^2, \mid e_j\in E^*,\;u_j\in U\} }.$$
\end{coro}

\begin{rema} Since $C_n/P_k=G_{\o}(k,2n)$ is a subvariety of the ordinary
  Grassmannian $G(k,2n)$, the $\PP^l$'s it contained are easy to
  describe: letting $M_0^a,N_0^b$ denote fixed linear spaces
  of dimensions $a$ and $b$, they are of the form
$$\{M_0^{k-1}\subset L\subset N_0^{k+l}\subset M_0^{\perp}\},$$
with $l\le 2n-2k+1$, or  
$$\{M_0^{k-l}\subset L\subset N_0^{k+1}\subset M_0^{\perp}\},$$
with $1<l\le k$ and $N_0$ isotropic. This second family of $\PP^l$'s
has two $C_n$-orbits, while the first family   breaks into a
 number of $C_n$-orbits that grows with $k$, 
indexed by the rank of the restriction of $\Omega$ to $N_0$. \end{rema}

\subsection{Odd spinor varieties}
 
We call  the homogeneous spaces $B_n/P_n$  the {\it odd spinor
varieties}. They are the usual $D_n$-spinor varieties seen as
$B_n$-homogeneous spaces. 

\begin{prop}
The tangent space and normal spaces of the odd spinor varieties 
$B_n/P_n$, as $H=A_{n-1}$-modules, are ($\tdim E=n$)\textsf{}:
$$\begin{array}{rcl}
T_1 &  =  & E^*, \qquad T_2 \;\;  =  \;\; \Lambda^2E^*, \\
N_p &  =  & \Lambda^{2p-1}E^*\oplus \Lambda^{2p}E^*. 
\end{array}$$ \end{prop}

\begin{coro}  
$$\tbase |\fff 2{B_n/P_n}|=
\ppp\overline { \{e\oplus e\we f,\;e\in E^*\backslash\{0\},\;f\in
E^*\backslash\CC  e\}}\simeq G(2,n).$$  
This base locus contains an open and dense $H$-orbit, the boundary of 
which is the union of the two (disjoint) closed $H$-orbits
 $$Y_1 = \ppp E^* =\PP(T_1)\quad and \quad Y_2\simeq G(2,E^*)
=G(2,n-1)\subset \PP(T_2).$$ 
The base locus of the $k$-th fundamental form is
$$\tbase |\fff k{B_n/P_n}|=
\ppp\overline{ \{e\oplus \Omega,\;e\in E^*\backslash\{0\},\;
{\rm rank}\, \Omega_{|e^{\perp}}<{\rm rank}\,\Omega\le 2k-2 \}}.$$  
\end{coro}

\proof Consider $V^{2n+1}\subset
W^{2n+2}$, and $E^n\subset F^{n+1}$, where $F$ is a null plane in $W$.
let $L=E\upperp\subset F$
Then $\La\bullet E^* \simeq \La{even}F^*$.
We may write $N_2=\La 3E^*\ot L
\oplus \La 4 E^*= \La 4 F= I_2(G(2,F))$. The analogous identities
hold for the higher normal spaces. After all, this is the
same projective variety as $D_n/P_n$.\qed

\begin{rema}\label{rbn} The varieties
parametrizing the $\PP^l$'s of $B_n/P_n=D_{n+1}/P_{n+1}$ are 
stratified as $B_n$-spaces. For $k\ge 2$, a connected component of
this variety is $D_{n+1}/P_{n,n-k-1}$.
If $\tilde V=\bcc{2n+2}$is endowed with a nondegenerate
quadratic form, a point of this space is a flag $F_0\subset M_0$
of isotropic subspaces of $\tilde V$ of respective dimensions $n-k-1$ and 
$n$. The associated line in $D_{n+1}/P_{n+1}$ is the space of
$n$-dimensional isotropic subspaces of $\tilde V$ containing $F_0$ and
cutting $M_0$ in dimension $n-1$. Let $V$ be the hyperplane of $\tilde 
V$ preserved by $B_n$. Then the $B_n$-orbits inside $D_{n+1}/P_{n,n-k-1}$
are indexed by the relative position of $F_0$ and $V$. There is a
closed orbit corresponding to $F_0\subset V$, isomorphic to
$B_n/P_{n-k-1}$, and its complement is an open orbit.  For $k=1$ or 
$k=3$, there is another connected component, parametrized by
$D_{n+1}/P_{n-1}$ and $D_{n+1}/P_{n-2}$ respectively. Each of them has
two $B_n$-orbits, the closed orbits being $B_n/P_{n-1}$ and $B_n/P_{n-2}$ 
respectively. \end{rema}
 
\section{Exceptional short roots and the octonions}

In this section we calculate $\tbase\ii$ for the exceptional
spaces corresponding to short roots, and give geometric
interpretations of these varieties and their linear spaces
in terms of the octonions. 

\subsection{$G_2/P_1$}

As an algebraic variety,   $G_2/P_1$  is
a familiar space, $G_2/P_1=\QQ^5\subset\pp 6$.
Studying it from an octonionic perspective will help us to
understand  
$F_4/P_4$   by analogy.

Identify  $\bcc 7\simeq Im \BO=V_{\lo 1}=V$.
Let $\phi\in\La 3 V^*$ be a generic element and
let $\rho :  GL (V)\ra GL(\La 3 V^*) $ be the induced
representation.  Here are some descriptions of 
  $G_2\subset GL(V)$ (see \cite{harvey} pp. 114, 116, 278 and 
  \cite{sv} chapter 2): 
$$\begin{array}{rcl}
G_2   &=&\text{Aut}(\BO)\\
&=&\{ g\in GL(V)  \mid  \rho (g)\phi =\phi\}  \\
&=& \{g=(g_+,g_-,g_0)\in \text{Spin}_8(V) \mid
g_+=g_-=g_0\}.\end{array}
$$
The third line should be understood as follows: let
$S_+,S_-,V_0$ denote the vector and two spin representations
of $Spin_8$ and choose appropriately an identification of
each of 
the three spaces with $\BO$, so that they are acted on by $g\in Spin_8$ in
three different ways, 
   call them  $(g_+,g_-,g_0)$. In this case (see \cite{harvey}
   p. 278), the {\em triality principle} of E. Cartan 
leads to the identification 
$$Spin_8= \{ (g_+,g_-,g_0)\in SO(\BO )\times SO(\BO )\times SO(\BO)
\mid g_+(uv)=g_-(u)g_0(v)\},$$ 
where $uv$ denotes octonionic multiplication.
When the three coincide one obtains an automorphism of the octonions,
showing the equivalence of  the second and third definition.
Harvey's description is explicit in bases. 
The connection between the first two interpretations is
that if one makes suitable identifications, for 
$u,v,w\in\; \text{Im}\OO$, we have
 $\phi (u,v,w)= Re [(uv)w]$.

The first definition is due to R. Bryant. It shows that $G_2$ is
not really an exceptional group, because it is defined by
a generic form. (Generic three forms on $\bcc m$ for $m>8$ are not
preserved by a positive dimensional group. For $m=6,8$, the groups
preserving such a form are classical.)

The third interpretation can be understood in terms of folding
Dynkin diagrams:

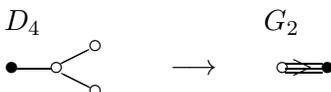
\setlength{\unitlength}{3mm}
\begin{picture}(15,5)(-15,-1.5)
\put(0,0){$\bullet$}
\put(2,0){$\circ$}
\put(3.7,1){$\circ$}
\put(3.7,-1){$\circ$}
\put(.45,.3){\line(1,0){1.65}}
\put(2.5,0.5){\line(2,1){1.25}}
\put(2.55,0.1){\line(2,-1){1.25}}
\put(0,2){$D_4$}
\put(7,0){$\lra$}
\multiput(12.45,.1)(0,0.2){3}{\line(1,0){1.65}}
\put(12,0){$\circ$}
\put(14,0){$\bullet$}
\put(12.7,0){$>$}
\put(11.5,2){$G_2$}
\end{picture}

\noindent This indicates that $G_2/P_1$ should be
  be understandable in terms of $D_4/P_1=\QQ^6$,
and in fact it is a generic hyperplane section.
 $\text{Im}\BO\subset \BO$ should be thought of as the
{\it traceless elements}, where the trace of an element is
its \lq\lq real\rq\rq\ 
 part and we call the hyperplane section $\{ \ttrace
=0\}$.
In what follows, $uv$ etc... refers to octonionic multiplication.

\begin{prop}\label{bl1}
Consider   $\QQ^5\simeq G_2/P_1=\ppp (Im\BO)_0\subset \PP (Im
  \OO)\simeq \PP (V_{\o_1})$.

Then $T_x G_2/P_1= T_1\op T_2\op T_3$
as an   $H= SL_2$-module.
Let  $A=\CC^2$, the standard representation of $SL_2$.
Then $T_1=A$, $T_2=\CC$ (the trivial representation) and $T_3=A^*$.
 Moreover, in a suitable normalization, 
$$\tbase |\FF^2_{G_2/P_1}|=  \PP \{a\,\op t\,\op a^* \mid
\langle a,a^*\rangle =  t^2 \}.$$
\end{prop}

\begin{prop}
 We have the following octonionic interpretations:
$$
\begin{array}{rcl}
G_2/P_1 &=&\{ [u]\in \ppp (Im \BO )\mid u^2=0\} \\
 \hat T_{[u]}G_2/P_1&=&\{ v\in Im \BO \mid uv+vu=0\}
=\{ v\in Im \BO \mid Re(uv)=0\} \\
\hat T_{[u]1}G_2/P_1 &=&\{ v\in Im \BO \mid uv =0\}.\ \end{array}
$$ 
\end{prop}

\begin{prop} 
The space $\BG (\pp 1, G_2/P_1)$
 of lines on $G_2/P_1$, has the following description:
$$
\BG (\pp 1, G_2/P_1)=\{ \ppp\{ u,v\}\mid 
 [u],[v]\in G_2/P_1\text{ such that }uv+vu= 0\}. 
$$
Note that $\BG (\pp 1, G_2/P_1)=G_o(2,7)$ and in particular
is of dimension seven.
\smallskip

The space  $G_2/P_2=\BG_{0}(\pp 1,\ppp ( Im\BO)_0)$
of $G_2$-homogeneous lines on
$G_2/P_1$ has dimension $5$, and admits  the following descriptions:

i. $G_2/P_2=\{ \ppp E\in \BG (\pp 1,\ppp V_{\lo 1}) \mid   E\lrcorner \phi=
0\}$

ii.   $G_2/P_2=\{  \ppp  E=\ppp\{ u,v\}
\mid    [u],[v]\in G_2/P_1, \   uv=0\}.$
\end{prop}

Proofs are left to the reader. The arguments are similar to,
but simpler than the arguments for the $F_4/P_4$ case below.
\smallskip

\begin{rema} Recall that $V_{\o_2}$ is the adjoint representation 
of $G_2$, so that $G_2/P_2$ is an adjoint variety. The relation 
with our description of $G_2/P_2$ as a space of special lines on the
quadric $G_2/P_1$ is as follows. Let $[u],[v]\in G_2/P_1$ be such that
$uv=vu=0$. One can then check that the map
$$d_{u,v}(z)=u(vz)-v(uz),\quad z\in\OO,$$
defines a nilpotent derivation of $\OO$,with $d_{u,v}^2=0$. \end{rema}

\medskip
Note that $G_{Q}(3,7)\simeq Q^6$, the space of $\pp 2$'s on $G_2/P_1$,
contains a special family of planes isomorphic to $G_2/P_1$ as follows:
  through each point of $G_2/P_1$  there is a plane
contained in $G_2/P_1$ tangent to $T_1$ (and it is completely
tangent to $T_1$ at this point only).  Perhaps it is better
to say the space of special $\pp 2$'s is parametrized by
$v_2(G_2/P_1)$ as the variety sits inside  $V_{2\o_1}\subset \La 3 V_{\o_1}$.
(Although $V_{\o_1}\subset\La 3 V_{\o_1}$,
the Veronse re-embedding is the correct factor as were there linear
spaces on the parameter space, they would determine larger linear spaces
on $G_2/P_1$ which do not exist.)

\subsection{The Cayley plane $\BO\pp 2$}

Let $\JO$ be the space of $3\times 3$ $\BO$-Hermitian symmetric matrices 
$$\JO =  \left\{ A =   \begin{pmatrix} r_1&\overline{x_3}&\overline{x_2} 
 \\
x_3&r_2&\overline{x_1} \\ x_2&x_1&r_3  \end{pmatrix},
r_i\in\CC, x_j\in\OO\right\}.
$$
$\JO$ can be equipped with the structure of  a {\em Jordan algebra} for
the commutative  product $A\circ B= \frac{1}{2} (AB+BA)$,
where $AB$ is the usual matrix product.
$\tdim_{\CC}\JO = 27$ and it is a model for the
$E_6$-module $V_{\lo 1}$.  There is a well-defined
{\em determinant} on $\JO$, which 
 is defined by same expression as the
classical determinant in terms of traces:
$$\det A =  \frac{1}{6}(\tr A)^3 -\frac{1}{2} (\tr A )(\tr A^2)
  +\frac{1}{3} \tr A^3.$$ 
  $E_6$ is  the
subgroup of $GL(\JO)= GL(27,\CC)$ preserving $det$.  
The notion of rank one matrices
is also well defined and the
{\it Cayley plane},  $E_6/P_1=\OO\PP^2 \subset \ppp (\JO )$ is the
projectivization of the rank one elements, with ideal the 
$2\times 2$ minors (see \cite{lansec}).

 Since $\a_1$ is not short, all linear spaces on $\BO\pp 2$ are described
by Tits geometries. In particular, $E_6/P_3$ is
the space of lines on $\BO\pp 2$ and $E_6/P_2$ is the space
of $\pp 5$'s on $\BO\pp 2$.

\subsection{$F_4/P_4=\BO\pp 2_0$} 

Here are some descriptions of $F_4\subset GL(\JO )$:
$$
\begin{array}{rcl}
F_4  & =&\{ g\in GL(\JO ) \mid  \ttrace ((\rho (g)A)^i) =\ttrace A^i 
\text{ for } i=1,2,3 \} 
\\
&=&\text{Aut}(\JO)
\\
&
 =& \{ g\in E_6  \mid
g_+=g_- \}\end{array}
$$
 The third
description is motivated by folding of
Dynkin diagrams:

 \setlength{\unitlength}{3mm}

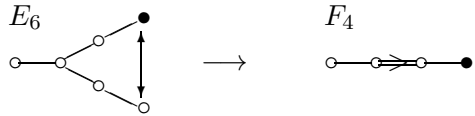
\begin{picture}(15,6)(-13,-2)

\multiput(0,0)(2,0){2}{$\circ$}
\put(0.45,.3){\line(1,0){1.65}}
\put(3.7,1){$\circ$}
\put(5.7,2){$\bullet$}
\multiput(3.7,-1)(2,-1){2}{$\circ$}
\put(2.5,.5){\line(2,1){1.3}}
\put(4.3,1.45){\line(2,1){1.4}}
\put(2.55,0.1){\line(2,-1){1.3}} 
\put(4.3,-0.9){\line(2,-1){1.4}} 
\put(0,2){$E_6$}
\put(8.2,0){$\lra$}
\multiput(14,0)(2,0){3}{$\circ$}
\put(20,0){$\bullet$}
\multiput(14.45,.3)(2,0){1}{\line(1,0){1.65}}
\multiput(18.45,.3)(2,0){1}{\line(1,0){1.65}} 
\multiput(16.4,.2)(0,0.2){2}{\line(1,0){1.7}} 
\put(16.8,0){$>$}
\put(14,2){$F_4$}
\put(5.9,.1){\vector(0,1){1.6}}
\put(5.9,.4){\vector(0,-1){1.6}}
\end{picture}\medskip

The equivalence  of the 
second and third descriptions can
 be proved by using the quadratic form
$\ttrace (A^2)$ to identify $\JO$ with $\JO^*$ and considering $g_+$ (resp.
$g_-$) as the two resulting elements of $GL(\JO)$. 
Harvey shows that the second definition implies the
first \cite{harvey} p. 296.   For the first definition,
one only needs two of the three forms to be preserved,
as any group preserving two preserves the third.

\smallskip
Geometric folding indicates $F_4/P_4$ should be understood
in terms of $\BO\pp 2$, and, as with $G_2/P_1$ above, it is
the hyperplane section $\{ \ttrace =0\}$. In what follows,
 $AB$ denotes the
usual matrix product of $A$ and $B$. Note that $A^2=A\circ A$.

\begin{prop}\label{bl2}
Consider   $\BO\pp 2_0= F_4/P_4\subset \PP (\JO_0)\simeq \PP (V_{\o_4})$.
Then  $T_x (F_4/P_4)= T_1\op T_2$ as an   $H= Spin_7$-module. 
Let  $U$ be the $7$-dimensional vector representation of $Spin_7$ and
$ {\mathcal S}(U)$ the spin representation, then $T_1= {\mathcal S} (U)$ and
$T_2=U$. The spinor variety $Y_1=\BS (U)$ is a six dimensional
quadric, and $Y_2=\QQ^5$ is a five dimensional quadric. Moreover,
we may identify $T_1\simeq\BO$ and $T_2\simeq Im\BO$ and with this
identification
$$
\begin{array}{rcl}
\tbase |\FF^2_{\BO\pp 2_0}| &=   \ppp
\overline { \{  (u,v) \in T_1\oplus T_2
\mid u\overline u=0, v\overline v=0, u\overline v=0  \} }
 &=\BS_5\cap H\end{array}
$$
where $\BS_5\cap H$ is a generic hyperplane section of the spinor variety
$\BS_5=D_5/P_5$. In particular, $\tbase |\FF^2_{\BO\pp 2_0}|$ is of 
dimension $9$, and is the closure of a $Spin_7$-orbit, the boundary 
of which is the disjoint union of $Y_1$ and $Y_2$. 
It is not homogeneous for any group.
\end{prop}

\proof 
The decomposition of the tangent space follows from \S 2.3.
Moreover, $\tbase |\FF^2_{\BO\pp 2_0}|$ must be a generic hyperplane section
of $\BS_5=\tbase |\FF^2_{\BO\pp 2 }|$
because $F_4/P_4$  is a generic hyperplane section of $E_6/P_6$. The explicit
description follows using the description of $T_x\BO\pp 2\simeq \BO\oplus\BO$
and the explicit description of $II_{\BO\pp 2}$ in \cite{lansec},
verifying that $\BO\op\BO_0$ is indeed a generic
hyperplane section or using the explicit description below.  Finally, we check that there is no
homogeneous space of dimension $9$ homogeneously embedded in a 
$\PP^{14}$.\qed

\begin{prop}\label{f4,2}
 We have the following octonionic interpretations:
$$
\begin{array}{rcl}
\BO\pp 2_0=F_4/P_4 &=&\{ [A]\in \ppp\JO_0\mid A^2=0\} \\
 \hat T_{[A]}\BO\pp 2_0&=&\{ B\in \JO_0 \mid A\circ B=0\} \\
\hat T_{1[A]}\BO\pp 2_0&=&\{ B\in \JO_0 \mid A  B=0\} \end{array}
$$
\end{prop}

\proof A calculation shows that an element
 $A\in\JO$ is rank one and traceless if and only if $A^2=0$.
Differentiation yields the second line.
 
To prove the third line, we 
first need to show that  if $[A]\in \BO\pp 2_0$
and $B\in \hat T_{[A]}\BO\pp 2_0$, the equation $AB=0$ is
$F_4$ invariant (although the matrix
product $AB$ is {\em not} $F_4$ invariant). Note that 
$F_4$ is generated by $SO_3$ and
$Spin_8$, where the action of $g\in SO_3$ is by
$A\mapsto gA{}^tg$, and that of $(g_+,g_-,g_0)\in Spin_8$ by
$$
\begin{pmatrix} r_1&\overline{x_3}&\overline{x_2} 
 \\
x_3&r_2&\overline{x_1} \\ x_2&x_1&r_3  \end{pmatrix}
\mapsto 
\begin{pmatrix} r_1&\overline{g_+(x_3)}&\overline{g_-(x_2)} 
 \\
g_+(x_3)&r_2&\overline{g_0(x_1)}
 \\ g_-(x_2)&g_0(x_1)&r_3  \end{pmatrix} .
$$
(This defines an automorphism of the Jordan algebra $\JO_0$ because
of the triality principle.) The $SO_3$ invariance is clear. 
Moreover, if we take
$$
A=\begin{pmatrix}i&1&0
 \\
1&-i&0 \\ 0&0&0  \end{pmatrix}
$$
then
$$
\hat T_{[A]}\BO\pp 2_0=
\begin{pmatrix} i\text{tr}(x_3) &\overline{x_3}& i\overline{x_1} 
 \\
x_3&-i\text{tr}(x_3)&\overline{x_1} \\ ix_1&x_1&0  \end{pmatrix}
$$
where $\text{tr}(u)=u-u_0=\frac 12(u+\overline u)$ is the \lq\lq
real\rq\rq\
part of $u$.
The $Spin_8$ invariance of the equation $AB=0$ is a straightforward
calculation, and follows again from the triality principle.  

With this model, $\{ B\in \hat T_{[A]}\BO\pp 2_0
\mid AB=0\}\simeq\{ x_1,\text{tr}(x_3)\}$ and we may consider
$\{x_1\}\subset \hat T_{[A]}/\{ A\}\simeq T$. Note that
$T$ is acted on by the subgroup of $Spin_8$ that preserves $A$,
which means that $g_+(1)=1$. By   \cite{harvey} p. 285,
$$
\begin{array}{rcl}Spin_7&=&\{ (g_+,g_-,g_0)\in
Spin_8\mid g_-=g_0\}\\
&=&
\{ (g_+,g_-,g_0)\in
Spin_8\mid g_+(1)=1\in\BO \}\end{array}
$$
(Note that this embedding of $Spin_7$ in $Spin_8$ is {\em not}
the standard one). 
Thus we explicitly see the $Spin_7=H$ action on $T$
and the decomposition of $T$ into $T_1\simeq \{x_1\}$
and $T_2\simeq \{ (x_3)_0\}$, respectively as the spin and
vector representations. In particular,  
  $\{ T_1+A\}=\{ B\in\hat T\mid AB=0\}$.\qed

\begin{prop}\label{f4p4,3}
The space $\BG(\pp 1,\BO\pp 2_0)$ of 
   lines on $\BO\pp 2_0$ has dimension $23$, and admits the following
description:
$$\BG(\pp 1,\BO\pp 2_0)=\{ \ppp\{  A, B\}\mid [A],[B]\in
\BO\pp 2_0\text{ such that }A\circ B= 0\}$$
 
The space  $F_4/P_3=\BG_{0}(\pp 1,\BO\pp 2_0)$ of
$F_4$-homogeneous lines on $\BO\pp 2_0$ has
dimension $20$, and admits the following   description:
$$
\BG_{0}(\pp 1,\BO\pp 2_0)= \{  \ppp\{ A, B\} \mid  
 [A],[B]\in \BO\pp 2_0\text{ such that }  AB=0\}.
$$
\end{prop}

\proof The geometric
descriptions of $\BG (\pp 1, \BO\pp 2_0)$ and
$\BG_{0} (\pp 1, \BO\pp 2_0)$ follow  immediately from
 proposition \ref{f4,2}, because 
$F_4/P_3$ is the space of lines on $F_4/P_4$ tangent to $Y_1$.
Moreover, 
$$\tdim (\BG(\pp 1,\BO\pp 2_0))=\tdim \BO\pp 2_0 +
\tdim (\tbase |\FF^2_{\BO\pp 2_0}|)-1= 15+9-1=23$$
verifies the dimension assertion.\qed

\medskip

\begin{prop} There are four   types of maximal (i.e. unextendable)
linear spaces through a point of $F_4/P_4$:

The space of $\pp 5$'s which is $5$-dimensional, parametrized
by the quadric $Q^5\subset\ppp T_2$.

A space of $\pp 4$'s which is $6$-dimensional, parametrized by
the quadric $Q^6\subset\ppp T_1$ or equivalently $\BS_{Q^5}$,
the variety of $\pp 2$'s in $Q^5\subset\ppp T_2$.

A space of $\pp 4$'s which is $6$-dimensional and having two components,
the two copies of $\BS_{Q^6}$, 
the variety of $\pp 3$'s in $Q^6\subset\ppp T_1$. 
\end{prop}

All other linear spaces can be deduced from these.

\proof
A $\pp k$  in $F_4/P_4$ corresponds to a $\pp{k-1}$ in $\BS_5\cap H$.
Let $L=\pp m\subset\BS_5\cap H$. The dimension $d_2$
of its projection $p_2(L)$ onto $Q^5\subset\ppp T_2$
is  $P_4$ invariant.

We choose a splitting $T=T_1\oplus T_2$ in order to use
the equations above describing $\tbase\ii$.
Relative to a choice of splitting,  $L$ is just the span of $p_1(L)$
and $p_2(L)$ so we can analyze $L$ accordingly.

Since $B_3$ acts transitively on $G_Q(k, T_2)$, we may choose
convienent $k$-planes to calculate with.

Without loss of generality, take
$v=\e_1+i\e_2\in T_2$, Then, writing $u=\a_0+\a_1\e_1+\hdots+\a_7\e_7$,
and using the standard octonionic multiplication table (e.g., see \cite{lansec})
the condition $uv=0$ implies
$u=i\a_3+\a_1\e_1+i\a_1\e_2 +\a_3\e_3+
i\a_5\e_4+\a_5\e_5 +\a_6\e_6 
+i\a_6\e_7$. In other words, we obtain a
$\pp 3_v\subset Q^6$ \lq\lq polar\rq\rq\ to $v$
which gives
rise to an unextendable $L^4=\langle v, \pp 3_v\rangle$.
 Note that we automatically have $u\overline u=0$.
Taking $M= \langle v,v'\rangle$ with
$v'= \e_6+i\e_7$, the additional condition $uv'=0$ implies
$u=\a  (\e_6 + i\e_7)$, i.e., is a point $q\in Q^6$.
A similar computation shows that any $q\in Q^6$ has
a $\pp 2$'s worth of points in $Q^5$ \lq polar\rq\ to it
 so we obtain
an unextendable $L^3=\langle \pp 2_q, q\rangle$.
If $p_2(L)$ is empty, then we are of course free to take
one of the two families of $\pp 3$'s on $Q^6$ as our maximal
linear space.\qed

\subsection{$F_4/P_3=\BG_{0}(\pp 1,\BO\pp 2_0)$} 
 
\begin{prop}\label{bl3}
  Consider $   F_4/P_3\subset \PP V_{\o_3}$.
Then   $T= T_1\oplus T_2\oplus T_3\oplus T_4$
as an   $H=     SL_3\times SL_2$-module. Let 
$\tdim E=   3$ and $\tdim U=   2$, then 
$$\begin{array}{c}
T_1=   E^*\ot U, \quad T_2=   E\ot S^2 U, \quad T_3=   U, 
\quad T_4=   E^*,\\ \hspace{2mm} \\
\tbase |\fff 2{\BG_{0}(\pp 1,\BO\pp 2_0) } |=    \ppp
\overline{ \{  e^{*}\ot u + e\ot u^2\in T_1\oplus T_2
\mid e^*\in E^*, e\in E, u\in U,\,  \langle e^*, e\rangle =   0  \} }.
\end{array}$$
This base locus $B$ is a nontrivial $\QQ^4$-bundle over $\ppp U=\pp 1$. 
In particular, $\tdim B=5$   and it has a dense open $SL_3\times SL_2$-orbit, 
the boundary of which  is the union of the two closed orbits 
$Y_1\subset\ppp T_1$ and $Y_2\subset\ppp T_2$. It is not homogeneous 
for any Lie group.\end{prop}
 
\begin{prop}\label{f4p3,2}
The space $\BG(\pp 1, F_4/P_3 )$ of 
$\pp 1$'s on $F_4/P_3$ has dimension $24$, and admits the following
description:
$$
\BG(\pp 1, F_4/P_3 )=\{ \{A\}\subset \{  A, B,C \}\mid [A],[B], [C]\in
\BO\pp 2_0\text{ such that } AB=AC=0,   B\circ C =0\}.
$$
The space $F_4/P_{2,4}=\BG_{0}(\pp 1 ,\BO\pp 2_0)$ of 
 $F_4$-homogeneous  $\pp 1$'s on $F_4/P_3$ has dimension $22$, and
admits the following description:
$$F_4/P_{2,4}=\BG_{0}(\pp 1,\BO\pp 2_0)=\{ \{A\}\subset \{  A, B,C
\}\mid [A],[B], [C]\in \BO\pp 2_0\text{ with }A  B= A  C=B  C =0\}.
$$
\end{prop}

\begin{coro} There are two types of maximal linear spaces passing through
a point of $F_4/P_3$;  the $\pp 3$'s corresponding to a $\pp 2$ in some
quadric in a fiber of
$\tbase |\fff 2{\BG_{0}(\pp 1,\BO\pp 2_0) } |$
considered as a fibration and the $\pp 2$'s corresponding to the $\pp 1$  
in the base.\end{coro}

\proof[Proof of \ref{f4p3,2}] We have $F_4/P_3\subset G(2,26)$ so
a line on $F_4/P_3$ must be a line of the Grassmanian as
well. Lines on $G(2,26)$ are determined by the choice
of a flag $\pp 0\subset\pp 2$. Here we need   $[A]= \pp 0\in F_4/P_4$
in both cases. 

In the first case $AB=AC=0$, $B\circ C=0$ are  necessary
and sufficient conditions that the line be contained in $F_4/P_3$,
as by \ref{f4p4,3}, we need $A(sB+tC)=0$ and
$(sB+tC)^2=0$ for all $[s,t]\in\pp 1$. Moreover, 
$\tdim \BG(\pp 1, F_4/P_3 )=24$
because the choice of $[A]$ is $15$ dimensions and then one
needs an element of $G_o(2,8)$, which
is of dimension $9$. (Here $\bcc 8\simeq T_1$.) 

In the second
case, considering $F_4/P_{2,4}$ as a $\pp 2$-bundle over $F_4/P_2$,
the conditions $AB=AC=BC=0$ follow from picking an element of
$F_4/P_2$, and the choice of $[A]$ is a choice of an element in the
fiber.\qed

\smallskip
\proof[Proof of \ref{bl3}] First, 
  $\tdim B=5$   because
$$\tdim (F_4/P_3) + \tdim \tbase |\fff 2{F_4/P_3} | -1=\tdim
\BG (\pp 1, F_4/P_3).$$ Moreover, we know that $B$ 
contains $Y_1$ and $Y_2$ and is irreducible.

Consider now  $[y_1+y_2]\in B$, 
with $y_j\in Y_j\subset\ppp T_j$. Write $y_1=   e\ot u$ and $y_2=   
e^*\ot v^2$. Conditions for such a point to belong to $B$ can 
only come from   components of 
$$T_1^*\ot T_2^*=   (E\ot E^*)\ot (U^*\ot S^2U^*) 
=   (\CC\op\fsl (E))\ot (U^*\op S^3U^*).$$ 

Suppose   that $\fsl (E)\ot S^3U^*$ were contained in $N_2^*$. Since 
$e\ot e^*$ is not a homothety, this would force $uv^2$ to be zero
in $S^3U$, hence $u$ or $v$ to be zero.  
If this   component were in $N_2^*$, then  $B$ would be included in 
$\PP T_1 \sqcup \PP T_2$, and would not be irreducible. 

Suppose now that $\fsl (E)\ot U^*$ were contained in $N_2^*$.
This set of equations would   force $u$ and $v$ to be parallel
because under 
the contraction $U\ot S^2U\ra U$,  $u\ot v^2$ maps to 
$\o (u,v)v$, where $\o\in\La 2U^*$. 
Similarly, the component $S^3U^*$ would force $\langle e,e^*\rangle=   0$. 

In conclusion, 
$\ppp\{ e\ot u\op e^*\ot u^2\mid \langle e,e^*\rangle=   0\}\subseteq B$.
Since both sets are irreducible of dimension five, the second one must be
the closure of the first one. 

The quadric bundle structure is given by
the application  $B\ra \ppp U$ defined by $[e^{*}\ot u + e\ot u^2]\mapsto
[u]$. This is a nontrivial bundle structure.   
Finally, to see that $B$ cannot be homogeneous,
note that there are no homogeneous nontrivial quadric fibrations
in dimension five. 
\qed

\medskip
Further calculations along this line show that each variety of
linear spaces is a finite union of $F_4$-orbits.

\vspace{2cm}

\begin{tabular}{lll}
Joseph M. Landsberg & \hspace*{1cm} & Laurent Manivel \\
School of Mathematics, & & Institut Fourier, UMR 5582 du CNRS \\
Georgia Institute of Technology, & & Universit\'e Grenoble I, BP 74 \\
Atlanta, GA 30332-0160 & & 38402 Saint Martin d'H\`eres cedex \\
USA & & FRANCE \\
{\rm E-mail}: jml@math.gatech.edu & &
{\rm E-mail}: Laurent.Manivel@ujf-grenoble.fr 
\end{tabular}
\end{document}